\documentclass[12pt]{article}

\usepackage{calrsfs}
\DeclareMathAlphabet{\pazocal}{OMS}{zplm}{m}{n}
\usepackage[letterpaper, nohead,dvips]{geometry}
\usepackage{amsmath}
\usepackage{amsfonts}
\usepackage{amssymb}
\usepackage{enumerate}
\usepackage{amsthm}
\usepackage{setspace}
\usepackage[all,cmtip]{xy}
\usepackage{graphicx}
\usepackage{chngcntr}
\usepackage[english]{babel}
\usepackage[section]{placeins}
\usepackage{float}
\usepackage{tikz-cd}
\usepackage{subfig}

\floatstyle{boxed}
\restylefloat{figure}

\renewcommand{\geq}{\geqslant}

\newcommand{\Rset}{\mathbb{R}}

\newcommand{\bL}{\pazocal{L}}
\newcommand{\bI}{\pazocal{I}}

\theoremstyle{plain}
\newtheorem{theorem}{Theorem}[section]
\newtheorem{proposition}[theorem]{Proposition}

\theoremstyle{definition}
\newtheorem{definition}[theorem]{Definition}
\newtheorem{remark}[theorem]{Remark}

\numberwithin{equation}{section}
\begin{document}

\title{In search of a new economic model determined  by logistic growth}

\author{Roman G. Smirnov\footnote{e-mail:  Roman.Smirnov@dal.ca}\,   and Kunpeng Wang\footnote{e-mail: kunpengwang@dal.ca} \\Department of
  Mathematics and Statistics\\
Dalhousie University\\ Halifax, Nova Scotia, Canada
  B3H~3J5}

\maketitle

\begin{abstract}
In this paper\footnote{Preliminary report, published in arXiv: https://arxiv.org/abs/1711.02625.} we extend   the work by Ryuzo Sato devoted to the development of economic growth models within the framework of the Lie group theory. We propose a new growth model based on the assumption of logistic growth in factors. It is employed to derive new production functions and introduce a new notion  of wage share. In the process   it is shown that the new  functions compare reasonably well against relevant economic data. The corresponding problem of maximization of profit under conditions of perfect competition is solved with the aid of one of these functions. In addition,  it is explained in reasonably rigorous mathematical terms why Bowley's law no longer holds true in post-1960 data. 
\end{abstract}

%%%%%%%%%%%%%%%%%%%%%%%%%%%%%%%%%%%%%%%%%%%%%%%%%%%%%%%%%%%%%%%%%%%%%%%%%%%%%%%
%\listoffigures

\section{Introduction}
\label{S1}

As is well known, a  production function is an essential feature of an economics growth model. Such a function can either be fixed, so that it is used to estimate the dynamics of other quantities, or it is the essential output of   the model, obtained by studying the dynamics of the input factors. 

%This project was originally conceived as a rather natural continuation of Ruyzo Sato's theory of technical change and economic invariance \cite{RS1981} in which we extended the assumptions about the Lie group properties of the technical progress representing the growth in factors of the production function in question. In our opinion, however,  the results obtained go beyond the scope of merely economic studies by means of mathematics (see Section \ref{S9}). Therefore, for brevity,  we dropped the adjective ``econoimc'' before ``growth'' in the original title of this paper. Nevertheless, we shall build our presentation around some important problemes that arise in economic theory. 

An example of the former application of a production function is the celebrated  Solow-Swan economic growth model \cite{RMS1956,TWS1956, G-PT2009}  introduced in the 1950s to  explain long-run economic growth, at which point  it  also  generalized and extended the Harrod-Domar model \cite{EDD1946, RFH1939} tasked with this undertaking prior.  The model in turn was later used as a starting point for the development of other economic growth models that emerged as its generalizations  (see, for example, Ferrara and Guerrini \cite{MFLG2009} and the relevant references therein). 

At the core of the Solow-Swan model and its generalizations is a production function $Y(t) = f(K(t), L(t))$, normally of the Cobb-Douglas type \cite{CWCPHD1928}, where the factors $K(t)$ and $L(t)$ represent capital and labor respectively. The function $Y(t)$ is required to satisfy the so-called  Inada conditions \cite{KII1963}. From a mathematical standpoint, the Solow-Swan economic growth model and its generalizations, for example, the   Ramsey-Cass-Koopmans model \cite{CD1965, TCK1965, FPR1928}, are governed  by a single nonlinear  differential equation or a system of such equations that describe the evolution of per capita  capital stock, consumption, etc.  

The theory of technical change and economic invariance developed by Ruzyo Sato \cite{RS1981} is an example of the latter approach, in which a production function is an output obtained within the framework of  a model. In particular, the author and his collaborators have derived the Cobb-Douglas production function as a consequence of the exponential growth in factors (capital and labor). 

 In this article we continue the development of Sato's theory  by  changing the assumptions about the Lie group theoretical properties of the technical progress representing the growth in factors. 

Recall that in 1928  Charles Cobb and  Paul Douglas published a paper \cite{CWCPHD1928} devoted to the study of the growth of the American economy during the period 1899-1922. To model the production  output they used the following function, introduced earlier by Knut Wicksell: 

\begin{equation}
Y  = A K^{\alpha}L^{\beta},
 \label{CD}
\end{equation}
where $K(t)$ and $L(t)$ are as before (i.e., in economic terms they are the factors of production), while $Y$ denotes  the total production, $A$ is  total factor productivity, and  $\alpha, \beta \ge 0$  are the output elasticities of capital and labor respectively. Sometimes the Cobb-Douglas function displays constant return to scale, which holds if 
\begin{equation}
\alpha+\beta=1, \, \alpha, \beta \ge 0.
\label{elasticity}
\end{equation}
The Cobb-Douglas function (\ref{CD}) can be easily derived under the assumptions that there is no production if either capital or labor vanishes, the marginal productivity of capital is proportional to the amount of production per unit of capital (i.e.,  $\frac{\partial Y}{\partial K} = \alpha \frac{Y}{K}$), and the marginal productivity of labor is proportional to the amount of production per unit of labor (i.e.,   $\frac{\partial Y}{\partial L} = \beta \frac{Y}{L}$).

More recently,  Ryuzo Sato \cite{RS1980, RS1981} (see also Sato and Ramachardan \cite{SR2014}),  while resolving the so-called Solow-Stigler controversy \cite{RMS1957, GS1961}, developed a Lie group theoretical framework to study technical progress and production functions. It can be viewed as an analogue of the Felix Klein approach to geometry formulated in his celebrated Erlangen Program \cite{FK1872} in which Lie transformation  groups play a central role. For instance, within this framework the Cobb-Douglas production function (\ref{CD}) can be recovered as an invariant of the one-parameter Lie group action \cite{AC1911}  that afford exponential growth in both $K$ and $L$  in the first quadrant of the two-dimensional Euclidean space $\Rset_{+}^2 = \{(K, L)| K, L \in \Rset_{+}\}$. The key idea employed by Sato \cite{RS1980, RS1981}, as well as Sato and Ramachadran \cite{SR2014} was to identify the corresponding exogeneous technical progress with the action of a one-parameter Lie group that acts in $C^2(\Rset_{+}^2)$.  
More specifically,  a Klein geometry can be described as  a pair $(G, H)$ where $G$ is a Lie group and $H$ is a closed Lie subgroup of $G$ such that the (left) coset space $G/H$ is connected. The group $G$ is called the principal group of the geometry and $G/H$ is called the space of the geometry, which is a homogeneous space for $G$.  For instance, in this view the pair $(SE(3), SO(3))$  describes the Euclidean geometry of $\Rset^3$ and  its objects, say,  surfaces can be classified modulo the action of the continuous isometry group $SE(3)$ (see, for example, Horwood {\em et al} \cite{HMS2005}, as well as Cochran {\em et al} \cite{CMS2017} for more details). By analogy,  a neoclassical growth model in the sense of Sato can be viewed as a pair $(G, \Rset^2_+)$, where the one-parameter Lie group of transformations $G$ acting in $C^2(\Rset_{+}^2)$ represents the technical progress in question. So far in the literature $G$ has been considered  to be either of a uniform (neutral) factor-augmenting type, that is $G: \overline{K} = e^{\alpha t}K$, $\overline{L} = e^{\alpha t}L$, for some  $\alpha \ge 0$, or representing a non-uniform, biased type, that is $G: \overline{K} = e^{\alpha t}K$, $\overline{L} = e^{\beta  t}L$ for some $\alpha, \beta  \ge 0$, $\alpha \not=\beta$. Therefore it is assumed in both cases that the economy grows exponentially (as per the corresponding growths in capital and labor), which was a reasonable assumption in the past based on the existing data at the time. It might  no longer be the case,  however, which may be attested to the fact, for example, that  the Cobb-Douglas function can no longer be used to describe adequatly  the growth of the American economy over a long run, including the recent decades ---  in the same way as it was done by Cobb and Douglas for the period 1899-1922 \cite{CWCPHD1928} (see Section \ref{S7} for more details). 

The main goal of this paper is to use the existing model to  develop a new mathematical paradigm that can be used to study the  current state of economy. Accordingly, in what follows we will modify the economic growth models described by Sato  within the framework of the Lie group theory according to the present economic realities \cite{PA2004}. More specifically,  we will replace in a neoclassical growth model in the sense of Sato $(G, \Rset^2_+)$ a group $G$ representing an exponential growth with another one-parameter Lie group that describes a {\em logistic growth}: 

$$G: \mbox{exponential growth} \rightarrow \mbox{logistic growth}. $$

This idea is currently being exploited and developed from different angles and in different directions quite extensively in the literature by economists and mathematicians alike (see, for example, \cite{AB2006, AB2007, BC2016, BCP2016,  CD2012, LG2010, LG2010a, LG2010b,  FG2008, FG2008a, FG2008b, FG2009, FG2009a}), which  is quite natural, given that the resources on our planet are  limited.

Therefore our first  task is to modify a basic  growth model $(G, \Rset^2_+)$ as described above and then, following Sato's approach, derive a new production function that may replace the Cobb-Douglas function (\ref{CD}) in  any models considered within the new paradigm of logistic growth, which is a reasonable further development, given that, for example, ``... the US economy is not well described by a Cobb-Douglas aggregate production function ...'' (see Antr\`{a}s \cite{PA2004} for more details). 

  Next, we will test the new production function derived purely by mathematical methods against a more up-to-date  data to verify its suitability for being part of any new mathematical models. Finally, we will reconsider several classical examples utilizing the properties of the Cobb-Douglas production funciton by replacing it with our new production function derived via the Lie group theoretical  approach developed by Sato and discuss the new results obtained under the assumption of logistic growth. 

This paper is organized as follows. In Section \ref{S2} we lay the groundwork for  the introduction of a new growth model and  derivation of new production functions. Specifically, we review the Lie group approach introduced in \cite{RS1981} and employ it to rederive the Cobb-Douglas function (\ref{CD}). In Section \ref{S3}  we depart from the growth model described by  Sato  based on exponential growth and introduce  instead a new one --- based on the assumption that factors grow logistically.  In Section \ref{S4} we derive a new production function (\ref{LPF1}) within the framework of the growth model (\ref{G1}) introduced in Section \ref{S3}. Section \ref{S5} is devoted to solving the problem of maximization of profit under condition of perfect competition, using the new production function (\ref{LPF1}). In Section \ref{S6} we explain, using mathematical reasonings and the results obtained in preceeding sections,  why Bowley's law \cite{ALB1900, ALB1937} no longer holds true in post-1960 data. In the process we also derive  another production function (\ref{productionf1}) and a new modified wage sare (\ref{ws1}).  In Section \ref{S7} we use statistical analysis to investigate  how estimations of the new production function (\ref{LPF1}) compare to economic data.  In Section \ref{S8} we make concluding remarks and summarize our findings.

\section{A Lie group approach to the study of holothetic production functions} 
\label{S2}

In this section we will briefly review the Lie group theoretical approach developed by Sato to study holothetic production functions and employ it to derive the Cobb-Douglas production function (\ref{CD}). Consider a  growth model $(G,  \Rset^2_+)$, where $G$ is a continuous one-parameter group of transformations (see \cite{ SR2014, RS1980, RS1981} for more details). In order to show that the increases in efficiency of inputs due to technical progress can be explained by economies of scale,  Sato  interpreted technical progress as the action of a one-parameter Lie group of transformations, for which the production function $Y = f(K, L)$ was an invariant. Under this arrangement the resulting transformation representing technical progress and generated by $G$,   indeed,  preserves the isoquant map, i.e., maps  one isoquant (or, in mathematical terms,  a level curve of $Y$) to another, that is, technical progress has the same effect as economies  of scale. 

More specifically, let capital and labour affected by technical progress and measured in the efficiency units, $\bar{K}$ and $\bar{L}$, be  given by
\begin{equation}
\bar{K} = \lambda_1 K, \quad \bar{L} = \lambda_2 L, 
\label{TP}
\end{equation}
where $\lambda_1$ and $\lambda_2$ represent the effect of the exogenous  technical progress. Following Sato and Ramachardan \cite{SR2014}, let us remark that  if $\lambda_1 = \lambda_2$ the change generated by technical progress is Hicks-neutral. If technical progress is factor  augmenting and biased, then $\lambda_1 \not= \lambda_2$. The functions $\lambda_i$, $i=1,2$ may depend on $t$ only, or they may be functions of $K/L$, which would imply that the rate of technical progress on different rays are different, but the rate is constant on each of them. They functions $\lambda_i$,  $i=1,2$ can also be functions of $K$, $L$ and $t$, which would entail that the rate of technical progress will also vary along a ray. In  what follows, we will also require that the technical progress functions $\lambda_i$, $i= 1,2$ represent the action of a one-parameter Lie group. 

Consider now the case when both $\lambda_i  = \lambda_i (t)$, $i=1,2$, moreover, $\lambda_1(t) = e^{\alpha t}$, $\lambda_2 (t)= e^{\beta t}$, $\alpha, \beta \geq 0$. Note, if $\alpha =\beta$ the change generated by such technical progress is Hick-neutral. Clearly, the corresponding transformations 
\begin{equation}
\bar{K} =  e^{\alpha t}K, \quad \bar{L} = e^{\beta t} L 
\label{TPe}
\end{equation}
form a continuous one-parameter Lie group, which follows from the fact, for example, that transformation (\ref{TPe}) determines the flow
\begin{equation}
\sigma(t, (K,L)) = \begin{bmatrix}e^{\alpha t} &0\\ 0 & e^{\beta t} \end{bmatrix}\begin{bmatrix}K \\ L \end{bmatrix}
\label{flow}
\end{equation}
generated by the following vector field 
\begin{equation}
U = \alpha K\frac{\partial}{\partial K} + \beta L\frac{\partial}{\partial L}, 
\label{VFU}
\end{equation}
which generates the Lie algebra of the one-parameter Lie group $G = \{g\,|\, g = \sigma_t, t \in \mathbb{R}\}$, where $\sigma_t: \mathbb{R}^2 \rightarrow \mathbb{R}^2$ is determined by  (\ref{flow}) for each fixed $t \in \mathbb{R}^2$. 

More generally, suppose a technical progress $T$ is defined  by the functions $\phi$ and $\psi$ such that 
\begin{equation}
T_t: \quad \bar{K} = \phi(K, L, t), \quad \bar{L} = \psi(K,L,t),
\label{TP1}
\end{equation}
where $t$ is the technical progress parameter and the functions $\phi$,  $\psi$ are analytic and functionally independent. Moreover, let us also suppose  the family of transformations $T_t$ (\ref{TP1}) forms a one-parameter Lie group $G$.  Recall, that Sato observed in  \cite{RS1981}  that in this case  a production function $f$ is  holothetic  under a continuous one-parameter Lie group  transformatoin (\ref{TP1}) iff 
\begin{equation}
Uf = \xi (K,L)\frac{\partial f}{\partial K} + \eta(K,L)\frac{\partial f}{\partial L} = H(f),
\label{Cov}
\end{equation}
where $\xi (K,L)  = \left(\frac{\partial \phi}{\partial K}\right)_{t_0 = 0}$, $\eta (K,L)  = \left(\frac{\partial \psi}{\partial L}\right)_{t_0 = 0}$. The condition of holotheticity is crucial from the economic standpoint, because it assures that the isoquant map  (i.e., the family of level curves of $f$) is invariant under the transformation (\ref{TP1}) representing the technical change, which means that under $T$ isoquants are mapped onto isoquants and the techinical change in this case is transformed into a scale effect.  

For example, if $\xi = \alpha K$ and $\eta = \beta L$  in (\ref{Cov}), $\alpha \not= \beta$, $\alpha, \beta >0$,  which means $\lambda_1 = e^{\alpha t}$, $\lambda_2 = e^{\beta t}$ in (\ref{TP}), $H(f) \not=0$,  it is a straigforward calculation, using the method of characteristic, that the general solution to the partial differential equation (\ref{Cov}) is given by \cite{RS1981} (see also \cite{RS1977}) 
\begin{equation}
Y = f\left[K^{1/\alpha}Q\left(\frac{L^{\alpha}}{K^{\beta}}\right)\right], 
\label{PF}
\end{equation}
where $Q(\cdot )$ is an arbitrary function. 

The converse problem was also considered by Sato. Specifically,  he established necessary and sufficient conditions for the existence of  a technical progress that affords holotheticity of a given production function (see Lemma 4 in \cite{RS1981} on p. 34). 

Now let us derive the Cobb-Douglas function (\ref{CD}) within the framework of the model $(G, \mathbb{R}^2_+)$, where the one-parameter Lie group of transformations $G$ determines the exponential growth (\ref{TPe}). Consider the partial differential equation  (\ref{Cov}) with the coefficients $\xi$ and $\eta$ determined by (\ref{TPe}) for $\bar{K} =  e^{a t}K$, $\bar{L} = e^{b t} L$, $a,b\geqslant  0$. Clearly, we can determine a particular production function (\ref{PF}) by specifynig the function $H(f)\not=0$ in (\ref{TPe}). Since $G$ in this case defines an exponential growth, it is natural to impose the corresponing condition on $H(f)$ --- so that it is also subject to an exponential growth. Indeed, let $H(f) = cf$, $c \geqslant 0$. Therefore  we have
\begin{equation}
Uf = a K\frac{\partial f}{\partial K} + b L\frac{\partial f}{\partial L} = cf,
\label{Cov1}
\end{equation}
or, alternatively, we can solve instead the following partial differential equation
\begin{equation}
X\varphi =a K \frac{\partial \varphi}{\partial K} + b L \frac{\partial\varphi }{\partial L} + cf\frac{\partial \varphi}{\partial f} = 0,
\label{Cov2}
\end{equation}
where $\varphi (K, L, f) = 0$, $\partial \varphi /\partial f \not\equiv 0$ is a solution to (\ref{Cov2}), while $f$ is a solution to (\ref{Cov1}) and an invariant. Solving the corresponding sysetm of ordinary differential equations
\begin{equation}
\frac{dK}{a K} = \frac{dL}{b L} = \frac{df}{cf},
\label{DE}
\end{equation} 
using the method of characteristics, yields the function (\ref{CD}), where $\alpha = \alpha (a,b,c), \beta = \beta(a,b,c)$. Unfortunately, the elasticity elements in this case do not attain economically meaningful values like (\ref{elasticity}). To overcome this problem Sato in \cite{RS1981} adjusted the model accodingly. Specifically, he introduces the notion of the simultaneous holothenticity, which implies that a production function is holothetic under more than one type of technical change simultaneously. Mathematically, it means that a production function is an  invariant of an integrable distribution of vector fields $\Delta$ \cite{AF2002} on $\mathbb{R}^2_+$, each representing a technical change as per the formula (\ref{Cov1}) (or, (\ref{Cov2})). More specifically, let us consider the following two vector fields, for which  a function $\varphi (K, L, f)$ is an invariant:  
\begin{eqnarray}
X_1\varphi =K \frac{\partial \varphi}{\partial K} +  L \frac{\partial\varphi }{\partial L} + f\frac{\partial \varphi}{\partial f}=0, \nonumber \\[0.3cm]
X_2\varphi =a K \frac{\partial \varphi}{\partial K} + b L \frac{\partial\varphi }{\partial L} + f\frac{\partial \varphi}{\partial f}=0.
\label{Hol} 
\end{eqnarray}
Clearly, the vector fields $X_1$, $X_2$ form a two-dimensional integrable distribution on $\mathbb{R}_+^2$: $[X_1, X_2] = \rho_1 X_1 + \rho_2X_2$, where $\rho_1 = \rho_2 = 0$. The corresponding total differential equation is given by  (see Chapter VII, Sato \cite{RS1981} for more details)
$$(fL - bfL)dK + (afK-fK)dL + (bKL - aKL)df = 0,$$
or, 
\begin{equation}
(1-b)\frac{dK}{K} + (a-1)\frac{dL}{L} + (b-a)\frac{df}{f} = 0. 
\label{TDE}
\end{equation}
 Integrating   (\ref{TDE}), we arrive at a Cobb-Douglas function of the form  (\ref{CD}), where the elasticity coefficients
$$\alpha = \frac{1-b}{a-b}, \quad \beta = \frac{a-1}{a-b}$$ satisfy the condition of constant return to scale  (\ref{elasticity}). 

\begin{remark}
\label{r1}
Note that, in principle, we could have used only one vector field generating   a partial differential equation of the type (\ref{Cov1}). However, the resulting Cobb-Douglas function would have had the   parameters satisfying the condition $\alpha\beta <0$  (see (\ref{CD})). The latter constraint on the parameters $\alpha$ and $\beta$ in (\ref{CD}) is incompatible with the economic growth theory  main postulates. We suppose that exactly  for this reason Sato \cite{RS1981} introduced the concept of simultaneous holotheticity. This arrangement, in particular, allows us to generate two-input  Cobb-Douglas  functions of the type (\ref{CD}) depending on a wide range of  parameters $\alpha$ and $\beta$, which we can, for instance,   make to satify the condition $\alpha + \beta = 1$,  so that the function (\ref{CD}) displays constant returns to scale as in the example above. 
\end{remark} 

These considerations lead to a very important conclusion, namely the Cobb-Douglas function,   derived within the framework of the growth model $(G,  \Rset^2_+)$, where the Lie group $G$ is determined by the exponential growth (\ref{TPe}), is precisely a manifestation of this exponential growth, or, more succinctly, we have 

$$\mbox{exponential growth} \Rightarrow \mbox{the Cobb-Douglas function,} $$
which means that the  Cobb-Douglas function (\ref{CD}) is a consequence of exponential growth representing technical change. 

\section{From exponential to logistic  growth models} 
\label{S3}

\begin{figure}
  \centering
  \subfloat[World Gold Reserves from 1845 to 2013, in metric tonnes (Wikipedia \cite{W1}).]{\label{oil}\includegraphics[width=0.5\textwidth]{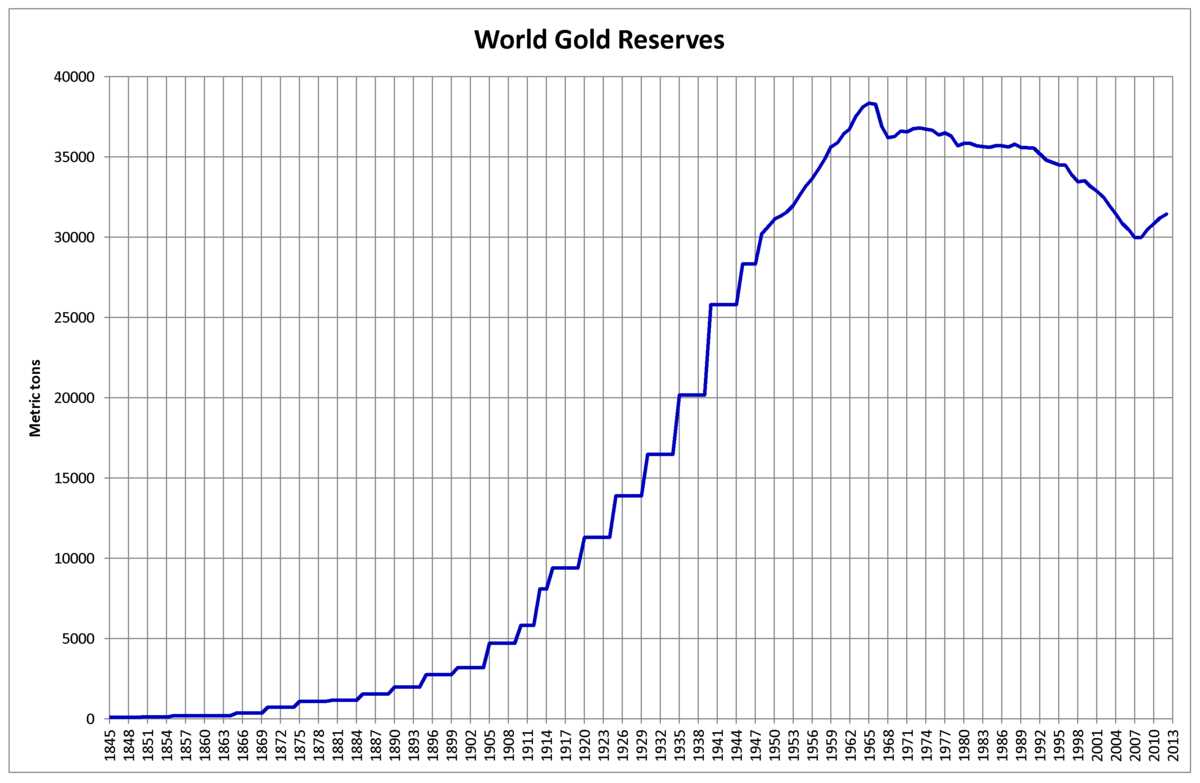}}
  \subfloat[World crude oil production 1930 to 2012 (Wikipedia \cite{W2}).]{\label{gold}\includegraphics[width=0.5\textwidth]{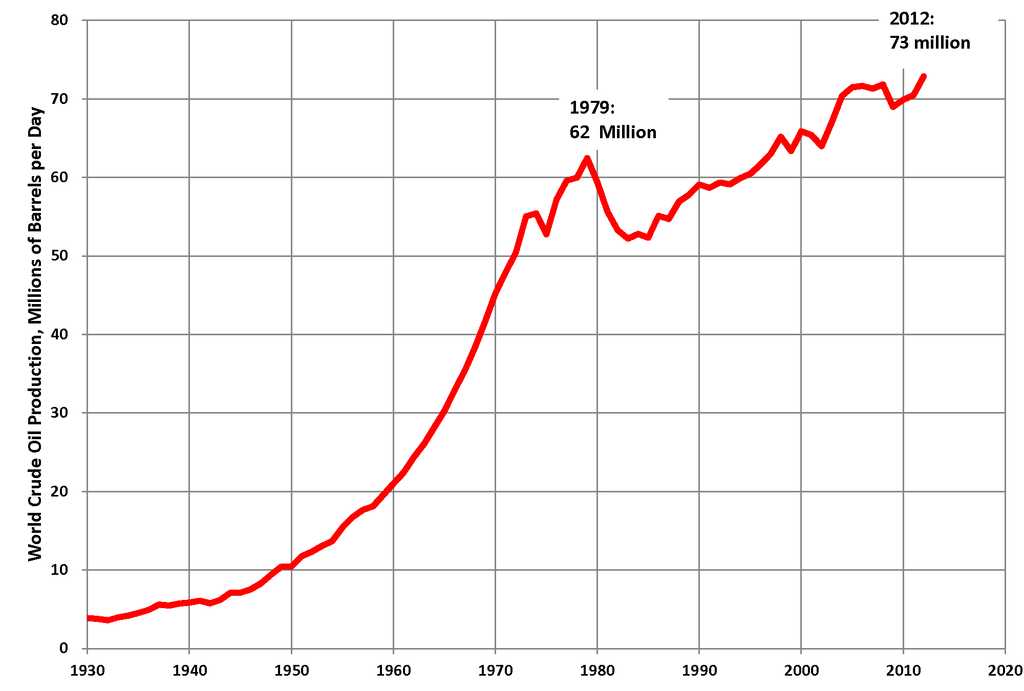}}
  \caption{\label{factors}Logistic growth in basic factors of production (gold and oil).}
\end{figure}

In this section we depart from the assumption that the input factors (i.e., capital and labor) grow exponentially in order to  extend  Sato's growth model $(G, \mathbb{R}^2_+)$. In what follows we assume that labor and capital grow logistically. There is already a substantial literature, starting with the pioneering paper by  Verhulst \cite{PFV1845}, in which the authors    have already based their considerations on this rather natural assumption, while studying various growth models with the aid of methods and techniques developed in  economics, mathematics and statistics (see, for example,  Brass \cite{WB1974},   Ferrara and Guerrini \cite{FG2008, FG2008a, FG2008b, FG2009, FG2009a, MFLG2009},  Leach \cite{DL1981}, Oliver \cite{ERO1982}, Tinter \cite{GT1952}).  The same assumption can be made about the growth in capital, if, for example, we look at such natural resources as oil  and gold as  proxies for energy and money respectively, it is quite evident that globally, given the fact that all resources are limited,  both the accumulation of gold reserves and oil production are subject to  logistic   rather than exponential growth, as can be illustrated by  Figure \ref{factors}.

We note that from the mathematical viewpoint it is also evident that there cannot be unbounded, continuous exponential growth, whether in terms of production, capital,  or population, on a planet with limited resources as per the following well-known theorem \cite{WR1976}: 
\begin{theorem}[Extreme value theorem]  If $K$  is a compact set and $f: K\to \mathbb{R}$ is a continuous function, then $f$  is bounded and there exist $p,q\in K$  such that $f(p)=\sup_{x\in K}f(x)$   and $f(q)=\inf_{x\in K}f(x)$.
\end{theorem}

In view of the above,  we propose the following  growth model based on the assumption that both capital  $K$  and labor $L$ are affected by logistic growth, namely 

\begin{equation}
(G_1, \mathbb{R}^2_+), \quad G_1: \bar{K} = \frac{N_KK}{K+\left(N_K-K\right)e^{-\alpha t}}, \quad  \bar{L} = \frac{N_LL}{L+\left(N_L-L\right)e^{-\beta t}}, 
\label{G1}
\end{equation}
where $\alpha, \beta >0$ and $N_K$, $N_L$ are the respective carrying capacities.  Clearly, $G_1$ is a one-parameter Lie group, acting in $\mathbb{R}_+^2$, whose flow is generated by the vector field 
\begin{equation}
U_1=  \alpha K\left(1-\frac{K}{N_K}\right)\frac{\partial }{\partial K} + \beta L \left(1-\frac{L}{N_L}\right) \frac{\partial }{\partial L}. 
\label{U1}
\end{equation}

\begin{remark}
It is also natural to  consider the  growth models $(G_2, \mathbb{R}^2_+)$ and  $(G_3,  \mathbb{R}^2_+)$  determined  by the assumption that    only  one of the two variables grow logistically, while the other is affected by exponential growth, that is
\begin{equation} 
(G_2, \mathbb{R}^2_+), \quad G_2: \bar{K} = \frac{N_KK}{K+\left(N_K-K\right)e^{-\alpha t}}, \quad  \bar{L} =e^{\beta t}L,
\label{G2}
\end{equation}\
or,
\begin{equation}
(G_3, \mathbb{R}^2_+), \quad G_3: \bar{K} = e^{\alpha t}K, \quad  \bar{L} = \frac{N_LL}{L+\left(N_L-L\right)e^{-\beta t}}.
\label{G3}
\end{equation} 
\end{remark} 
Following the approach developed by Sato in \cite{RS1981}, we can now determine the corresponding family of production functions by solving the partial differential equation determined by the vector field $U_1$ (\ref{U1}): 
\begin{equation}
U_1f =  \alpha K\left(1-\frac{K}{N_K}\right)\frac{\partial f}{\partial K} + \beta L \left(1-\frac{L}{N_L}\right)\frac{\partial f}{\partial L} = H(f),
\label{Covlogistic}
\end{equation}
where $H(f)$ is an arbitrary function of $f$. Employing the method of characteristics, we arrive at the following family of functions: 
\begin{equation}
Y = f_1\left\{\left(\frac{K}{\left|N_K-K\right|}\right)^{1/\alpha}Q\left[\left(\frac{L}{\left|N_L- L\right|}\right)^{\alpha}\left(\frac{\left|N_K-K\right|}{K}\right)^{\beta}\right]\right\},
\label{Y1}
\end{equation} 
where $Q(\cdot )$ is an arbitrary function. We note that for $N_K = N_L = 1$ and $K, L \ll 1$  the family of functions given by (\ref{Y1}) $f_1 \sim f$, where $f$ is given by  (\ref{PF}).  Therefore we arrive at the following
\begin{proposition}
The most general family of production functions holothetic within the  growth model (\ref{G1}) is given by (\ref{Y1}). 
\end{proposition} 

\begin{remark}
The same argument applied to the ``partially" logistic neoclassical growth models (\ref{G2}) and (\ref{G3}) yields the families of functions
\begin{equation}
\label{Y2}
Y = f_2\left\{\left(\frac{K}{\left|N_K-K\right|}\right)^{1/\alpha}Q\left[L^{\alpha}\left(\frac{\left|N_K-K\right|}{K}\right)^{\beta}\right]\right\}
\end{equation}
and
\begin{equation}
\label{Y3}
Y = f_3\left\{K^{1/\alpha}Q\left[\left(\frac{L}{\left|N_L- L\right|}\right)^{\alpha}{K}^{-\beta}\right]\right\},
\end{equation}
respectively. 
\end{remark}

Our next goal is to derive a new production function under the assumption of logistic growth in both capital $K$ and labor $L$. Since the Cobb-Douglas function (\ref{CD}) has been shown above to be a member of the family of production functions (\ref{PF}) determined within the neoclassical growth model $(G, \mathbb{R}^2_+)$, where the Lie group $G$ is given by  (\ref{TPe}), it is natural to seek a new production function compatible with the logistic growth determined by the action of the Lie group $G_1$ (\ref{G1}) within the growth model $(G_1, \mathbb{R}^2_+)$. This is the subject of the considerations that follow. 

\section{From logistic growth to a new production function} 
\label{S4}

In  Section \ref{S2} we saw how the Cobb-Douglas production function could be derived   as an element of the family of production functions (\ref{PF})  within the framework of the growth model $(G, \mathbb{R}^2_+)$, where the Lie group $G$ was  defined  by (\ref{TPe}). Now let us consider the new growth model $(G_1, \mathbb{R}_+)$, where the Lie group $G_1$ was  given by (\ref{G1}). Before we formally derive the corresponding production function as an element of the family of production functions (\ref{Y1}), following the procedure outlined above, let us first give a reasonable justification for the calculations that we shall present below. 

Recall that a   necoclassical growth model of the Solow type may be  defined as follows (see, for example,  Jones and Scrimgeour \cite{CIJDS2008}, a model with decay in produced capital was studied in Cheviakov and Hartwick \cite{CH2009}): 

\begin{eqnarray}
Y & = & f(K, L; t), \nonumber \\
Y & = & C+I , \nonumber \\
I & = &sY, s >0\\
\dot{K} & = & I - \delta K, \, K_0, \, \delta \ge 0, \nonumber \\
L & =&  L_0e^{\alpha t}, \, L>0, \, \alpha \ge 0, \nonumber 
\label{model}
\end{eqnarray} 
where $C$ and $I$ represent consumption and investment (savings) respectively, while $\delta$ denotes  depreciation of capital. It is also assumed that the production function $f$ satisfies the Inada conditions \cite{KII1963}: 
\begin{enumerate}
\item $f_K, f_L >0$, this condition accounts for growth in both $K$ and $L$, 

 \item $f_{KK}, f_{LL} <0$,  that implies diminishing marginal returns also in both $K$ and $L$, 

\item $f$ has constant returns to scale, that is $f(\lambda K, \lambda L) = \lambda f (K, L)$ for all $\lambda >0$, 

\item $f$ satisfies  the following properties: 

$$\lim_{K\to 0}f_K = \infty, \lim_{K\to\infty} f_K = 0,  $$  
$$\lim_{L\to 0}f_L = \infty, \lim_{L\to\infty} f_L =0.$$
\end{enumerate}
For example, the Cobb-Douglas function (\ref{CD}) satisfies the above assumptions, provided the condition (\ref{elasticity}) holds. Such a model and its generalizations ensure steady long-run growth, ignoring short-run fluctuations. Since the pioneering  paper by Solow \cite{RMS1956}  was published in 1956 the model (\ref{model}) and its many generalizations  have played the most prominent role in  the development of the  endogenous growth theory. Clearly, the production function $Y$ is the cornerstone of the model and if it satisfies the Inada conditions the growth is driven by decreasing marginal returns from the very beginning for all $K, L>0$.  Many important examples of endogenous growth support this assumption (see, for example,  Cobb and Douglas \cite{CWCPHD1928}). Nevertheless, there are situations when  growth cannot be described by a strictly concave production function. For instance, at a microeconomic level a company may develop a product based on an original idea, such a product initially can be sold unrestricted in the absence of competition, generating increasing marginal returns. After a while, a competition may become a factor  (e.g., other companies may introduce similar products)  affecting the sales of the original product, whose market share may shrink. In turn, this situation in a long-run will manifest itself in  decreasing marginal returns.  Mathematically, the corresponding  production function will no longer be strictly concave. Capasso {\em et al} \cite{CED2012} gave a different motivation for the introduction of  a  (globally) nonconcave production function based on the idea of ``poverty traps". The authors also pointed out  two examples of models based on nonconcave production functions: Skiba \cite{AKS1978} (economics) and Clark \cite{CWC1971} (mathematical biology).  A macroeconomic example of such a scenario of growth can be found in Tainter \cite{JAT1988} (see Figure   16, p. 109). 

To address the issue Capasso {\em et al} \cite{CED2012} (see also Engbers {\em et al} \cite{EBC2014},  La Torre {\em et al} \cite{LLM2015}, Anita {\em et al} \cite{ACKL2013, ACKL2015, ACKL2017}) employed  a purely heuristic approch to  introduce a new general family of    production functions of the form
\begin{equation}
\label{LPF}
Y = f_4(K, L) = \frac{\alpha_1K^pL^{1-p}}{1 + \alpha_2K^pL^{1-p}}
\end{equation}
reducible to the Cobb-Douglas function (\ref{CD}) and enjoying an ``$S$-shaped'' (concave-convex) behavior for $p\ge 2$.  Clearly, the functions of the class (\ref{LPF}) have a horizontal asymptote as $(K, L) \to (\infty, \infty)$ when $\alpha_2 \not=0$  and  are compatible with logistic growth. These functions were used by the authors as a cornerstone for building a new, highly non-trivial generalization of the Solow model with spacial component in which they did not make assumptions about logistic growth for $L$. It is worth mentioning at this point that  Ferrara and Guerrini \cite{FG2008, FG2008a, FG2008b, FG2009, FG2009a, MFLG2009}, while generalizing the Ramsey and Solow models of economic growth, assumed logistic growth in $L$, but kept the Cobb-Douglas function (\ref{CD}) intact. 

The introduction of the family of production functions (\ref{LPF}) is certainly a big step in the right direction, nevertheless these functions cannot account for all possible examples of growth (and decay). For example, a production function can exhibit growth, followed by a period of stabilization  and then decay (see, for example, \cite{JBC1973}). Another option is growth, followed by a period of stabilization, which is followed by growth again. In this view our next goal is to derive a more general production function that can be used to describe a wider range of economic growth models, including the situations outlined above. We shall employ the  Lie group theoretical method  developed by Sato \cite{RS1981} and briefly described in Section \ref{S2}. 

Indeed, consider the growth model $(G_1, \mathbb{R}_+^2)$ given by (\ref{G1}). Next, we are going to identify a member of the family (\ref{Y1}) compatible with logistic growth given by (\ref{G1})  by imposing the corresponding constraints on the RHS of the equation (\ref{Covlogistic}). By analogy with the case of the Cobb-Douglas function derived by Sato \cite{RS1981} within the framework of the growth model $(G, \mathbb{R}_+^2)$, where the action of the Lie group $G$ is determined by (\ref{TPe}), let us consider the following partial differential equation determined by the vector field  $U_1$ given by  (\ref{U1}): 
\begin{equation}
U_1f =  a K\left(1-\frac{K}{N_K}\right)\frac{\partial f}{\partial K} + b L \left(1-\frac{L}{N_L}\right)\frac{\partial f}{\partial L} = c f\left(1-\frac{f}{N_f}\right),
\label{Covlogistic1}
\end{equation} 
or, in other words, let us specify the function $H(f)$ in (\ref{Covlogistic}) to be  $ c f\left(1-\frac{f}{N_f}\right)$ that implies logistic growth in the production function as well. Compare (\ref{Covlogistic1}) with the equation (\ref{Cov1}). 
\begin{remark}
We note that the choice for the RHS of (\ref{Covlogistic1}) is  not arbitrary. It turns out that in order to obtain a meaningful solution one needs to assure that the properties of the function $H(f)$  in (\ref{U1}) are compatible with the logistic growth determined by (\ref{G1}). For example, if we set  $H(f) = f$ in (\ref{U1}), which would imply that the growth in both $K$ and $L$ is logistic, while $f$ grows exponentially, the resulting production function would have singularities (see  the equation (\ref{productionf2})). Therefore the above equation reflects the fact  that the growth determined by (\ref{Covlogistic1}) is consistent for all quantities involved, that is for $K$, $L$ and $f$. 
\end{remark} 

Next, we employ the same reasoning that Sato in \cite{RS1981} based his derivation of the Cobb-Douglas function (\ref{CD}) upon (see also Section \ref{S2}). Let us assume that the production functions  in two sectors of an economy (or, two countries) are identical, so that the aggregate production function sought is of the same form. However, it does not necessarily mean that the technical changes in both sectors are also  the same. That is in what follows we shall give conditions under which the aggregate production function in question  is holothetic  under two types of technical changes simultaneously and solve (again) the corresponding {\em simultaneous holotheticity problem}. In mathematical terms, let us consider the following two vector fields acting on  a function $\varphi (K, L, f)$:  
\begin{eqnarray}
X_3\varphi & = &K \left(1-\frac{K}{N_K}\right)\frac{\partial \varphi}{\partial K} +   L \left(1-\frac{L}{N_L}\right)\frac{\partial\varphi }{\partial L} + f\left(1-\frac{f}{N_f}\right) \frac{\partial \varphi}{\partial f}=0, \nonumber \\
X_4\varphi & = &aK \left(1-\frac{K}{N_K}\right)\frac{\partial \varphi}{\partial K} +  b L \left(1-\frac{L}{N_L}\right)\frac{\partial\varphi }{\partial L} + cf\left(1-\frac{f}{N_f}\right) \frac{\partial \varphi}{\partial f}=0.\nonumber
\end{eqnarray}
Clearly, the vector fields $X_3$ and $X_4$ form an integrable distribution $\Delta$ on $\mathbb{R}_+^2$, because $[X_3, X_4] = \rho_3X_3 + \rho_4X_4$, where $\rho_3 = \rho_4 = 0$. Then the corresponding total differential equation which has $\varphi (K, L, f) = \mbox{const}$ for a solution assumes the following form: 
$$
\begin{array}{rcl}
\displaystyle \left[(c-b)f\left(1-\frac{f}{N_f}\right)L\left(1 - \frac{L}{N_L}\right) \right]dK   &+ & \\[0.5cm]
 \displaystyle \left[(a-c)f\left(1-\frac{f}{N_f}\right)K\left(1 - \frac{K}{N_K}\right) \right]dL &+ &\\[0.5cm]
\displaystyle \left[(b-a)f\left(1-\frac{K}{N_K}\right)L\left(1 - \frac{L}{N_L}\right) \right]df  &= & 0,
\end{array} 
$$
or, 
\begin{equation}
(c-b)\frac{dK}{K\left(1- \frac{K}{N_K}\right)} + (a-c)\frac{dL}{L\left(1- \frac{L}{N_L}\right)} + (b-a)\frac{df}{f\left(1- \frac{f}{N_f}\right)}=0.
\label{TDE1}
\end{equation}

Integrating the differential equation (\ref{TDE1}) (compare it with (\ref{TDE})), we arrive at a solution of the form  $\varphi (K, L, f) = 0$ defined in the open domain 
$$D = ]0, N_K[\times]0,N_L[\times ]0,N_f[ \subset \mathbb{R}^3$$ and satisfying the condition $\frac{\partial \varphi}{\partial f} \not\equiv 0$. Solving for $f$ by the impilcit function theorem, we arrive at the following hypersurface in $\mathbb{R}^3$: 
\begin{equation}
\label{LPF1}
 Y = f_{5}(K, L) = \frac{N_{f_5}K^{\alpha}L^{\beta}}{C\left|N_K- K\right|^{\alpha}\left|N_L - L\right|^{\beta} + K^{\alpha}L^{\beta}}, \, (K,L) \in \mathbb{R}_+^2,
\end{equation}
where $C \in \mathbb{R}$ is the constant of integration, $\alpha =\frac{c-b}{a-b} $, $\beta = \frac{a-c}{a-b} $. Note $\alpha + \beta = 1$. Note that in view of the symmetry of the differential equation (\ref{TDE1}), we  could have solved the equation $\varphi (K, L, f) = 0$ for $K$ and $L$ as well. The function $Y = f_5(K, L)$ given by (\ref{LPF1}) whose range is $]0, N_{f}[$  coinsides with the function  $\varphi (K, L, f) = 0$ on $D$. 

%extending
Furthermore, we note that in the subset $D' = ]0, N_K[\times ]0, N_L[ \subset \mathbb{R}_+^2$  of the domain of the function $Y = f_5(K, L)$ its growth is governed  by the logistic growth in the factors $K$ and $L$. Note that in this region the growth of the production function $f_5$ is ``$S$-shaped'',  which agrees with the assumptions that led to the introduction of the production function (\ref{LPF}).  However, the production function (\ref{LPF1}) is also  defined outside of the region $D'$, which impies in turn that its shape in the subset $\mathbb{R}_+^2 \setminus D' = [N_K, \infty[\times[N_L, \infty[$ is determined by the growth in $K$ and $L$ that goes beyond the respective carrying capacities $N_K$ and $N_L$. We will elaborate on this matter  without loss of generality while dealing with the corresponding one-input analog of the new two-input production function (\ref{LPF1}) below.

%, that is we have 
%$$\dot{K} = - K\left(1- \frac{K}{N_K}\right), \quad \dot{L} = - L\left(1- \frac{L}{N_L}\right),$$
%where $K \ge N_K$, $L \ge N_L$.

 %In other words, for  $K \ge N_K$, $L \ge N_L$ the growth in $K$ and $L$ respectively is of the  ``overproduction'', ``increased participation",  or ``overshoot''  type (see Section \ref{S9}). Note that the production function $f_5$ in the region  $]N_K, \infty[\times]N_L, \infty[$ actually decreases, attaining its absolute  maximum at $(K, L) = (N_K, N_L)$ (see Remark \ref{r12}). More specifically, the condition $K \ge N_K$ may mean the infusion of ``cheap money" (loans), that in turn increases participation (say, more people able to property than afforded by economic conditions), hence the condition $L \ge N_L$. 
%extending

We conclude, therefore, that by analogy with the algorithm  based on the Lie group theory methods  devised  by Sato and  applied in \cite{RS1981} to generate the Cobb-Douglas function (\ref{CD}), we have used it, after some modifications, to generate a {\em new production function (\ref{LPF1})}. More succinctly, we have
$$\mbox{logistic growth} \Rightarrow \mbox{the new production function (\ref{LPF1}).} $$ 
\begin{remark}
\label{r12}
%We note first that  the new production function $f_5$ given by (\ref{LPF1}) attains  its  absolute  maximum at $(K, L) = (N_K, N_L)$: $f_{5} (N_K, N_L) = N_f$. The next observation is of particular importance.

 Taking the limit as $K, L \to \infty$ (even though $K$ and $L$ cannot grow beyond a certain ``horizon" - see below), we obtain
\begin{eqnarray}
\lim_{\substack{K\to \infty\\ L\to \infty}}f_5(K, L) &= & \lim_{\substack{K\to \infty\\ L\to \infty}}\frac{N_{f_5}K^{\alpha}L^{\beta}}{C\left|N_K- K\right|^{\alpha}\left|N_L - L\right|^{\beta} + K^{\alpha}L^{\beta}}\\
                                                                                 &=& \lim_{\substack{K\to \infty\\ L\to \infty}}\frac{N_{f_5}}{C\left|\frac{N_K}{K}-1\right|^{\alpha}\left|\frac{N_L}{L}-1\right|^{\beta}+ 1}\\
                                                                                  &=& \frac{N_{f_5}}{C+1}.
\end{eqnarray} 
The quantity
\begin{equation}
\label{ss}
S_{f_5} = \frac{N_{f_5}}{C+1}
\end{equation}
 is the {\em steady state} of the new production function $f_5$ given by (\ref{LPF1}). Note that by changing the constant $C$ in (\ref{ss}) we can regulate the steady state $S_{f_5}$. 

\end{remark}
\begin{remark}
\label{r2}
See Remark \ref{r1}. 
\end{remark}

\begin{remark}
\label{r3}
We observe that the new production function $f_5$ (\ref{LPF1}) is reducible to the production function (\ref{LPF}) proposed by  Capasso {\em et al} \cite{CED2012} when $K$ and $L$ $\ll$ $N_K$ and $N_L$ respectively, $N_L, N_K \approx 1$, $C = 1$ in (\ref{LPF1}) and $ \alpha_1 = N_{f_5}$,   $ \alpha_2 = 1$ in (\ref{LPF}) .
\end{remark}

\begin{remark} 
\label{r4}
Figure \ref{LPF1surface} presents the surface  of a two-input production function of the type (\ref{LPF1}) for $N_f = 120$, $\alpha  = \beta = 3$, $N_K = 113$, $N_L = 115$, $C = 1.18$ without singularities (see Remark \ref{r5}). 
\begin{figure}
  \caption{A two-input  production function of the type (\ref{LPF1}) with isoquants.}
\label{LPF1surface}
  \centering
   \includegraphics[width=\textwidth]{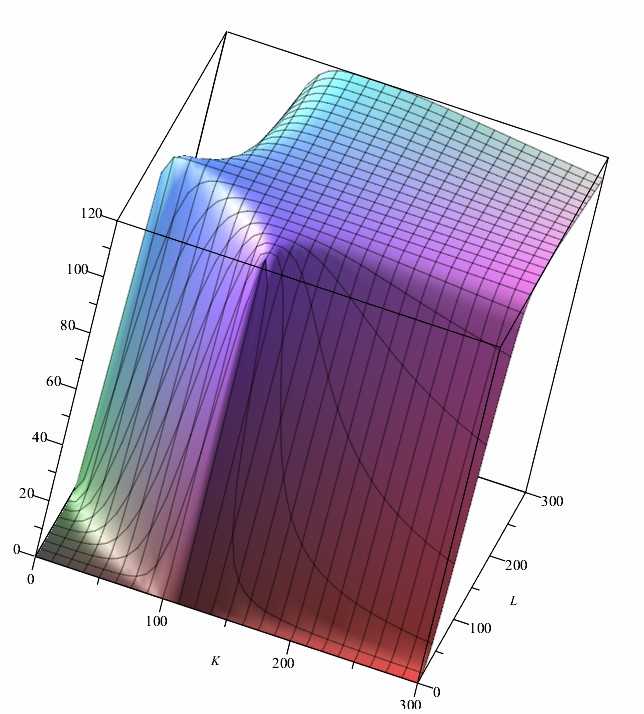}
\end{figure}

\end{remark} 

\begin{remark}
\label{r5}
Employing  the same procedure,   we  can determine now in a fairly straightforward  manner the corresponding one-input analogue  of  the new two-input production function (\ref{LPF1}).  Thus, let us derive a new production function $Y = f(x)$ whose growth is governed the growth in   $x$ which we assume to be logistic. Hence, we can formulate the following problem within the framework of the growth model $(G_2, \mathbb{R}_+)$:  
\begin{equation}
(G_2, \mathbb{R}_+), \quad G_2:  \bar{x} =  \frac{N_xx}{x+\left(N_x-x\right)e^{-a t}}, \, a >0, x \in \mathbb{R}_+, 
\label{fl}
\end{equation}
\begin{equation}
U_2 f =  a x\left(1-\frac{x}{N_x}\right)\frac{d f}{d x}  = b f\left(1-\frac{f}{N_f}\right), 
\label{Covariant3}
\end{equation} 
where the vector field $U_2 = ax\left(1-\frac{x}{N_x}\right)\frac{\partial }{\partial x}$ represents the infinitesimal action defined by the Lie group $G_2$ (\ref{fl}). Separating the variables and integrating the differential equation (\ref{Covariant3}) yields the follwoing solution (production function): 
\begin{equation}
\label{newfunction1}
Y = f_{6} (x) =\frac{N_{f_6}x^{\alpha}}{C|N_x-x|^{\alpha} + x^{\alpha}}, 
\end{equation}
where $C \in \mathbb{R}$ is the constant of integration and $\alpha = b/a$ with the corresponding steady state given by 
\begin{equation}
S_{f_6} = \frac{N_{f_6}}{C+1}.
\label{Sf6}
\end{equation} 
%extending
Note that in this case as well the new production function (\ref{newfunction1}) exhibits first an ``$S$-shaped'' growth in the region $]0, N_x[$, followed by a decline for $x > N_x$. Let us investigate this case from the economics point of view in more detail.
%extending

Let us recover the corresponding group action that affects the input $x(t)$, so that this action could be viewed as growth which entails the condition $\dot{x}(t) >0$. Indeed, consider the infinitesimal action $\tilde{U}$ given by $\tilde{U}=\tilde{U}_1 \frac{\partial}{\partial x}+\tilde{U}_2 \frac{\partial}{\partial y}$ so that $\tilde{U}f_6 = 0.$ Solving the last equation, we obtain the following solutions: 
\begin{equation}\label{eq3}
\begin{array}{ll}
U_1&=a \dfrac{x(N_x-x)}{N_x}, \\[0.3cm]
U_2&=b \dfrac{y(N_y-y)}{N_y}
\end{array}
\end{equation}
and
\begin{equation}\label{eq4}
\begin{array}{ll}
U_1&=a \dfrac{x(x-N_x)}{N_x},\\[0.3cm]
U_2&=b \dfrac{y(y-N_y)}{N_y}.
\end{array}
\end{equation}
In view of the fact that $x(t)$, $y(t)$ $>0$, it follows from \eqref{eq3} and \eqref{eq4} that
\begin{equation}\label{eq18}
\begin{array}{rcl}
\dot{x}&=a \dfrac{x(N_x-x)}{N_x}, &{0<x<N_x},\\[0.3cm]
\dot{y}&=b \dfrac{y(N_y-y)}{N_y},  & {0<y<N_y}
\end{array}
\end{equation}
and
\begin{equation}\label{eq9}
\begin{array}{rcl}
\dot{x}&= a \dfrac{x(x-N_x)}{N_x}, & {x>N_x}, \\[0.3cm]
\dot{y}& = b \dfrac{y(y-N_y)}{N_y}, & {y>N_y}, 
\end{array}
\end{equation}
so that both $x(t)$ and $y(t)$ represent growth. Solving the above equations, we obtain
\begin{equation}\label{eq5}
x(t)=\left\{
\begin{array}{rcl}
&\dfrac{N_x}{1+C_1e^{-a t}},  & 0<x(t)<N_x,\\[0.3cm]
&\dfrac{N_x}{1+C_2e^{a t}},  & x(t)>N_x,
\end{array}\right.
\end{equation}
where $C_1>0$ and $C_2>0$ are constants of integration. Next, we determine the time interval corresponding to growth in $x(t)$. It follows (\ref{eq5})  that $t>0$ for $0<\frac{N_x}{1+C_1e^{-a t}} <N_x$ and $0<t<\frac{1}{a}\ln \frac{1}{C_2}$ for $  \frac{N_x}{1+C_2e^{at}}   >N_x$. Substituting the equation (\ref{eq5}) into (\ref{newfunction1}), we arrive at the following function: 
\begin{equation}\label{eq6}
y (t) =\left\{
\begin{array}{lrc}
&\dfrac{N_{f_6}}{C(C_1e^{-a t})^{\alpha}+1},  & 0<t<t_1,\\[0.3cm]
&\dfrac{N_{f_6}}{C(C_2 e^{a t})^{\alpha}+1}, & t_1<t<\frac{1}{a} \ln\frac{1}{C_2},
\end{array} \right. 
\end{equation}
where $t_1$ is the time at which the function shifts from the logistic to a different growth type. Let us assume $\alpha$ to be a positive integer. Furthermore, we note that $y(t)$ increases or decreases depending on whether $\alpha$ is odd or even respectively.  To assure that (\ref{eq6}) is compatible with (\ref{newfunction1}) we assume that $\alpha$ is an even integer (see below). Next, rewrite the production function given by (\ref{eq6}) as follows:
\begin{equation}
\label{Heaviside}
y=(H_0(t)-H_{t_1}(t))y_1(t)+H_{t_1} y_2(t),
\end{equation}
where $H_c(t)$ is the Heaviside (unit) step function, $$y_1(t) = \frac{N_{f_6}}{C(C_1e^{-a t})^{\alpha}+1}, \quad y_2(t) = \frac{N_{f_6}}{C(C_2 e^{a t})^{\alpha}+1}.$$
In this view the function (\ref{Heaviside}) may interpreted  as an impulse response function. Indeed, a sudden change in the input at $t = t_1$ causes a jump in the output  from $y_1(t)$ to $y_2(t)$. From the economics viewpoint we can identify this phenomenon as a ``shock" \cite{CAS1980}, which means that a sudden change in exogenous factors yields the corresponding  sudden change in production (see \cite{KPP1996, PS1998, AHJ2014} for more details and referenses). The gap between $y_1(t)$ and $y_2(t)$ caused by a sudden change in $x(t)$ at $t=t_1$ is given by 
\begin{equation}
\begin{array}{ll}
d_{(y_1,y_2)}(t_1)= \dfrac{C N_{f_6}(C_2^a e^{b t_1}-C_1^a e^{-b t_1})}{(C(C_1e^{-a t_1})^a+1)(C(C_2 e^{a t_1})^a+1)},
\end{array}
\end{equation}
where $d_{(y_1,y_2)}(t_1)$ denotes the distance between the two curves at $t=t_1$. Next, we note that 
\begin{equation}
\label{Sy}
y (t) \to \frac{N_{f_6}}{C +1}, \quad \mbox{as} \quad   t \to \frac{1}{a} \ln \frac{1}{C_2}. 
\end{equation}
Note that if $\alpha$ is an even number, the RHS of (\ref{Sy}) is precisely the steady state (\ref{Sf6}).

 Figure \ref{LPF1curve}  presents the graph of a one-input production function of the type (\ref{newfunction1}) generated for $N_{f_6} = 100$, $\alpha = 2$ and $C = 2$. Note the function given by (\ref{newfunction1}) defines  an invariant $I (K, L)$ of the infinitesimal action determined by vector field $U_1$ (\ref{U1}) for $f_6 = K$ (or, $L$) and $x= L$ (or, $K$), namely $U_1 I = 0$, where 
$$I(K, L) = \frac{L^{\alpha}}{|N_L-L|^{\alpha}} \cdot \frac{N_K-K}{K}.$$ 
\begin{figure}
\caption{A one-input  production function of the type (\ref{newfunction1}).}
\label{LPF1curve}
\centering
\includegraphics[width=\textwidth]{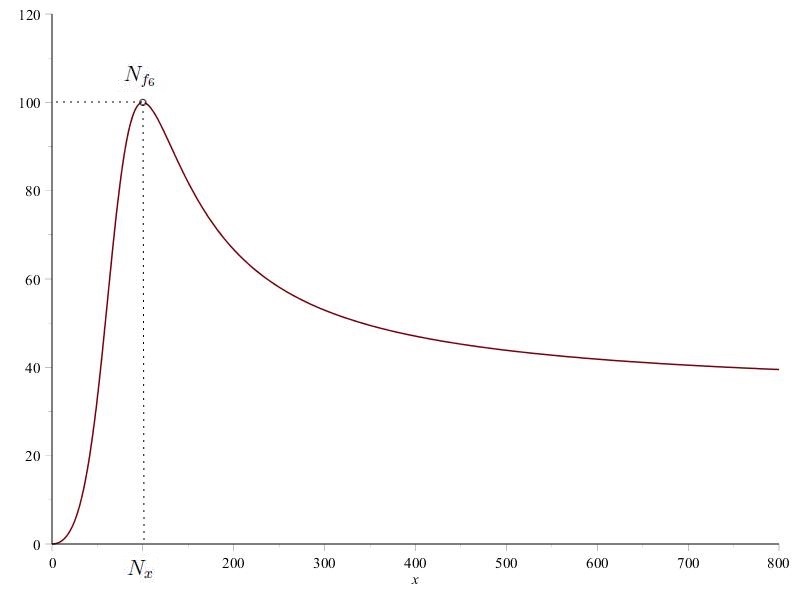}
\end{figure} 
\end{remark} 

\begin{remark}
\label{r6}
Repeating the above calculation within the frameworks of the growth models (\ref{G2}) and (\ref{G3}), we arrive at the production functions

\begin{equation}
Y = f_7(K, L) =  \frac{N_{f_7}K^{\alpha}L^{\beta}}{C\left|N_K- K\right|^{\alpha} + K^{\alpha}L^{\beta}}
\label{klogistic}
\end{equation}
and
\begin{equation}
Y=f_8(K, L) =  \frac{N_{f_8}K^{\alpha}L^{\beta}}{C\left|N_L - L\right|^{\beta} + K^{\alpha}L^{\beta}},
\label{llogistic} 
\end{equation}
respectively, where the parameters $\alpha$ and $\beta$ are the same as in (\ref{LPF1}). 

We also note that the functions (\ref{klogistic}) and (\ref{llogistic}) are elements of the families (\ref{Y2}) and (\ref{Y3}) respectively, as expected. 
\end{remark}

\section{The problem of maximization of profit under conditions of perfect competition}
\label{S5}

In 1947 Paul Duglas gave his presidential address  to the American Economics Association in which he referred to a  coherent assembly of the statistical  evidence accumulated  in the course of the previous 20 years while he and other people were  studying various economic data that confirmed the validity of the Cobb-Douglas production  function. It is safe to assume that this event marked the beginning of its  universal acceptance by the mainstream economic science.  He wrote  in \cite{PHD1976}: ``... the Cobb-Douglas function was being widely used, and that a host of younger scholars led by my former student, Paul Samuelson, his colleague Solow and Marc Nerlove, the son of my friend and former colleague, Samuel Nerlove, were all pushing forward into new and more sophisticated fields.''  In fact, Marc Nerlove gave  a series of lectures at the Econometric Workshop held at the University of Minnesota in 1957, which were subsequently published a few years later in  a book \cite{MN1965}. One of the problem considered by the author was the problem of maximization of profit of a firm under conditions of perfect competition in both factors and product markets under the assumption that the revenue of the firm from sales was determined by the Cobb-Douglas production function. In what follows we shall  solve  the problem using the same arguments {\em mutatis mutandis} as in \cite{MN1965} by assuming that the revenue of the firm from sales is now determined by the new production function (\ref{LPF1}). 

Consider  an individual firm functioning under  conditions of perfect competition in both factors and product markets. It attempts to maximize its profits by employing optimal quantities of inputs and producing an optimal quantity of output. At the same time its purchases of factors and supply of output do not affect the prices of the factors involved and the final product. Therefore  the said prices are assumed to be given, while the profits are to be maximized. Let $\Pi$, $p_0$, $p_1$, $p_2$ be the profit, the price of the final product, the cost of using one unit of capital, and  the wage of labor respectively. Hence, we have 
\begin{equation}
\label{profit} 
\Pi = p_0Y - p_1K - p_2L.
\end{equation} 
Traditionally, in problems like this the output $Y$ is assumed to be related to the inputs $K$ (capital) and $L$ (labor)  by the Cobb-Douglas production function (\ref{CD}). Instead, suppose now $Y$ is related to $K$ and $L$ via the new production function $f_5$ (\ref{LPF1}). Next, let us solve the problem of maximization of the profit $\Pi$ given by (\ref{profit}) subject to the constraint implied  by (\ref{LPF1}). The corresponding Lagrangian function $\cal L$ is readily found to be 
\begin{equation}
\label{lagrangian}
{\cal L}(Y, K, L, \lambda) = \Pi - \lambda \left(Y - \frac{N_{f_5}K^{\alpha}L^{\beta}}{C\left|N_K- K\right|^{\alpha}\left|N_L - L\right|^{\beta} + K^{\alpha}L^{\beta}}\right),
\end{equation}
where $\lambda$ is a Lagrange multiplier. For profit to be a maximum, the total differential 
\begin{equation}
\label{differential1}
d{\cal L}(Y, K, L, \lambda) = d(\Pi - \lambda g) = 0,
\end{equation}
where \begin{equation}\label{g}g = Y  - \frac{N_{f_5}K^{\alpha}L^{\beta}}{C\left|N_K- K\right|^{\alpha}\left|N_L - L\right|^{\beta} + K^{\alpha}L^{\beta}}.\end{equation} The condition (\ref{differential1}) yields
\begin{equation}\label{eq2}
\begin{array}{rcl}
\dfrac{\partial {\cal L}}{\partial \lambda}  &=&-Y+\dfrac{N_{f_5}K^{\alpha}L^{\beta}}{C\left|N_K- K\right|^{\alpha}\left|N_L - L\right|^{\beta} + K^{\alpha}L^{\beta}}=0,\\ [0.5cm]
\dfrac{\partial {\cal L}}{\partial K}&=&-p_1+p_0\dfrac{\beta N_{f_5} C (N_K-K)^{\alpha} K^{\alpha} (L^{\beta-1}(N_L-L)^{\beta}+(N_L-L)^{\beta-1}L^{\beta})}{(C(N_K-K)^{\alpha}(N_L-L)^{\beta}+K^\alpha L^\beta)^2}=0,\\[0.5cm]
\dfrac{\partial {\cal L}}{\partial L}&=&-p_2+p_0\dfrac{\alpha N_{f_5} C (N_L-L)^{\beta}L^{\beta}( K^{\alpha-1}(N_K-K)^{\alpha}+(N_K-K)^{\alpha-1} K^{\alpha})}{(C(N_K-K)^{\alpha}(N_L-L)^{\beta}+K^\alpha L^\beta)^2}=0,\\[0.5cm]
\dfrac{\partial {\cal L}}{\partial Y} &=&p_0-\lambda=0.
\end{array}
\end{equation}
The equations (\ref{eq2}) give us necessary conditions for maximum profit. Solving (\ref{eq2})  with the aid of the computer algebra system Maple, we get 
\begin{equation}
\label{solution}
\begin{array}{rcl}
Y &=& \dfrac{N_{f_5}K^{\alpha}L^{\beta}}{C\left|N_K- K\right|^{\alpha}\left|N_L - L\right|^{\beta} + K^{\alpha}L^{\beta}}, \\[0.5cm]
\alpha &=& \displaystyle \dfrac{p_2N_{f_5} K(N_K-K)}{p_0 N_K Y (N_{f_5}-Y)},\\[0.5cm]
\beta &=&\displaystyle \dfrac{p_0 N_K Y(N_{f_5}-Y)\ln\dfrac{|N_{f_5}-Y|}{CY}-p_2N_{f_5} K(N_K-K)\ln\dfrac{|N_K-K|}{K}}{p_0 N_K Y(N_f-Y)\ln\dfrac{|N_L-L|}{ L}}.
\end{array}
\end{equation}
The resulting equations (\ref{solution}) are sufficient to determine the variables $Y$, $K$ and $L$. The corresponding sufficient conditions for maximum profit are provided by the necessary conditions established above supplemented by the following second-order condition: 
$$d^2{\cal L} <0,$$
or, given the fact that $\Pi$ in (\ref{lagrangian}) is linear in $Y$,  $K$ and $L$ (see (\ref{profit})) and $\lambda = p_0$ by (\ref{eq2}), we have
\begin{equation}
\label{secondorder}
d^2\tilde{g} >0,
\end{equation}
where $$\tilde{g}(K,L) = \frac{p_0N_{f_5}K^{\alpha}L^{\beta}}{C\left|N_K- K\right|^{\alpha}\left|N_L - L\right|^{\beta} + K^{\alpha}L^{\beta}}.$$
Solving (\ref{secondorder}), using Maple, we arrive at the following set of inequalities: 
\begin{equation}
\begin{array}{l}
\alpha(\alpha-1) <0, \\
\beta(\beta-1)<0, \\
(2K-N_K)(2L-N_L)+N_L(2K-N_K)\beta+N_K(2L-N_L)\alpha  >0, \\
(2K-N_K)(2L-N_L)-N_L(2K-N_K)\beta-N_K(2L-N_L)\alpha  >0. 
\end{array}
\end{equation}
The first two inequalities entail that $0 < \alpha, \beta <1$. The second two inequalities imply that   $K>N_K/2$ and $L>N_L/2$. Hence, we arrive at the following conditions that assure maximum profit: 
\begin{equation}
\label{max2}
\begin{array}{l}
0 < \alpha, \beta <1, \quad K > N_K/2,  \quad L>N_L/2, \\
(2K-N_K)(2L-N_L)+N_L(2K-N_K)\beta+N_K(2L-N_L)\alpha  >0, \\
(2K-N_K)(2L-N_L)-N_L(2K-N_K)\beta-N_K(2L-N_L)\alpha  >0. 
\end{array}
\end{equation}
Next, we observe that since  $\lim_{t\to\infty}K(t) = N_K$ and $\lim_{t\to\infty}L(t) = N_L$, the last inequality in  (\ref{max2}) implies that
\begin{equation}
\label{last}
0< \alpha+\beta <1,
\end{equation} 
which in turn implies that the assumption of perfect competition and maximization of profit are inconsistent  in the case when 
$$\alpha + \beta \ge 1.$$

Finally, we conclude that the equations and inequalities  (\ref{solution}), (\ref{max2}) and (\ref{last}) constitute   sufficient conditions for maximum profit of a firm in the environment of perfect competition. The equations (\ref{solution}) determine the output a firm will deliver and the inputs of factors it will employ once the prices of the product and factors are established. Therefore the conclusions are pretty much the same as in the case when the revenue is determined by the Cobb-Douglas production function (\ref{CD}) considered in Nerlove \cite{MN1965}. The case of imperfect competition in both factor and production markets will be considered in a forthcoming paper. 

 Note that all of the calculations above have been carried out under the assumption that $C>0$. If $C<0$ the condition (\ref{last}) changes to $\alpha +\beta >1$. 

\section{The wage share and logistic growth} 
\label{S6}

The labor share is the fraction of national income, or the income  of a particular economic sector,  defined as the share  which is payed out to employees. Therefore it is often also called the wage share.   As is well-known, the  wage  share in the economic growth models governed by the Cobb-Douglas production  function (\ref{CD}) is a constant. More specifically, its constant value can be derived directly from the Cobb-Douglas function and  expressed in terms of the output elasticity of capital  in a simple and elegant way when the Cobb-Douglas function, say,  enjoys constant return to scale  (see, for example, Rabbani  \cite{SR2017}). The invariance of the wage share is subject to Bowley's law \cite{ALB1900, ALB1937} or the law of the constant wage share, which states that the share of national income  that is paid out to the employees as compensation for their work (normally, in the form of wages), remains unchanged (invariant) over time \cite{JMK1939, HMK2011, DS2011}. Economic data collected in different countries till about 1980 gave rise to and most strongly supported this law, which was widely accepted by the economics community at the time. However, this is no longer the case on both counts (see, for example, Schneider \cite{DS2011} for more details and references).

In view of the mathematics presented  above, it should not be viewed as a surprise. Indeed, the ivariance of wage share is linked to the Cobb-Douglas production function, which in turn is a consequence of exponential growth, as shown by Sato \cite{RS1981}. Next, since one of the the main points of this research project is the idea that we must depart from the exponential growth model and accept the logistic one, let us ivestigate how this transition affects the wage share. 

In what follows we shall  propose a new formula for the wage share compatible with logistic growth and support our claim by a rigorous mathematical analysis. 

First, let  us recover the formula for the wage share as an invariant of a prolonged 
infinitesimal group action given in terms of the corresponding projective coordinates defined as the output-capital ration $Y/K = y$ and the labor-capital output $L/K  = x$. The terminology and notations that we will use  are compatible with those adopted by Olver \cite{PJO1993, PJO1995} and Saunders \cite{DJS1989}. Consider a general production function 
\begin{equation}
\label{Yg}
Y = f (K, L; t)
\end{equation}
under the assumption that the dependent and independent  variables $K$, $L$ and $Y$ grow exponentially: 
\begin{equation}
\label{growth}
\bar{K}=Ke^{\alpha t}, \quad \bar{L}=Le^{\beta t}, \quad \bar{Y}=Y e^{\epsilon t}, \,\, \alpha, \beta, \epsilon \geqslant 0. 
\end{equation}
In view of the material presented in Section \ref{S2} we know that the production function (\ref{Yg}) is bound to be of the Cobb-Douglas type (\ref{CD}), in terms of the projective coordinates it assumes the following form: 
\begin{equation}
Y = f(x; t).
\label{Yp}
\end{equation} 
Clearly, the one-parameter  Lie group  of transformations (\ref{growth}) induces the corresponding action on the projective coordinates, which is also exponential: 
\begin{equation}
\bar{y}=ye^{\gamma t}, \quad \bar{x}=xe^{\lambda t}, \,\, \gamma, \lambda \geqslant  0
\end{equation}
with the corresponding infinitesimal action given by the  vector field $\mathbf{u}$ (compare it with (\ref{VFU})) given by
\begin{equation}
\mathbf{u}=\lambda x \frac{\partial}{\partial x}+\gamma y \frac{\partial}{\partial y}.
\end{equation}

Following Saunders \cite{DJS1989}, let us suppose that  $(\mathbb{R}^2,\pi,\mathbb{R})$ is a trivial bundle so that $\pi=pr_1$ and $(x,y)$ are   adapted coordinates. Then the corresponding jet bundles are $(J^1 \pi, \pi_1,  \mathbb{R})$ and $(J^1\pi, \pi_{1,0}, \mathbb{R}^2),$ as per the  commutative diagram (\ref{diagram}), where the first-jet manifold of $\pi$ is given by
\begin{equation}
J^1 \pi=\left\{j^1_p \phi: p \in \mathbb{R}, \phi \in \Gamma_p(\pi) \right\}
\end{equation}
with  adapted coordinates $(x,y,y_x).$
\begin{equation}
\label{diagram}
\begin{tikzcd} [column sep=5em, row sep=5em]
J^1 \pi \arrow{r}{\pi_{1,0}} \arrow[swap]{d}{\pi_1} & \mathbb{R}^2 \arrow{d}{\pi} \\
\mathbb{R} \arrow{r}{id} & \mathbb{R}
\end{tikzcd}
\end{equation}
Here $\pi_1=\pi \circ \pi_{1,0}.$

Next, the first prolongation of $\mathbf{u}$ on $\mathbb{R}^2$ is the following vector field  $\mbox{pr}^{(1)}\mathbf{u}={\bf u}^{(1)}$, which  has to be    a symmetry of the Cartan distribution on $J^1 \pi$ (see  Saunders \cite{DJS1989} for more details),  that is the vector field
\begin{equation}
\mbox{pr}^{(1)}\mathbf{u} =  {\bf u}^{(1)} = \lambda x \frac{\partial}{\partial x}+\gamma y \frac{\partial}{\partial y}+\xi(x,y,y_x) \frac{\partial}{\partial y_x}
\end{equation}
 is required to be a symmetry of the Cartan distribution on $J^1 \pi$. Indeed, consider  a basic contact form $\omega = dy - y_x dx.$  Next, in view of the above, we require that the one-form ${\bL}_{\mathbf{u}^{(1)}}(\omega)$ is a contact form, where $\bL$ denotes the Lie derivative. Thus, we compute 
\begin{equation}
\label{prol1}
\begin{array}{ll}
{\bL}_{\mathbf{u}^{(1)}} (\omega) & = {\bL}_{\mathbf{u}^{(1)}} (dy-y_xdx) \\
&= {\bL}_{\mathbf{u}^{(1)}}(dy)-({\bL}_{\mathbf{u}^{(1)}} y_x)dx-y_x ({\bL}_{\mathbf{u}^{(1)}}(dx))\\
&=\mbox{d}(\mathbf{u}^{(1)})(y))-(\mathbf{u}^{(1)} (y_x))\mbox{d}x-y_x\mbox{d}(\mathbf{u}^{(1)}(x)) \\
&= \gamma \mbox{d}y-\xi(x,y,y_x)\mbox{d}x-\lambda y_x\mbox{d}x \\
&=\gamma(\omega+y_x\mbox{d}x)-\xi(x,y,y_x)-\lambda  y_x \mbox{d}x \\
&= \gamma \omega+(\gamma y_x-\xi(x,y,y_x)-\lambda y_x)\mbox{d}x.
\end{array}
\end{equation}
The last line of (\ref{prol1})  implies that the expression in the parentheses above vanishes, which entails that $\xi(x,y,y_x)=(\gamma-\lambda)y_x$. Therefore the first prolongation $\mathbf{u}^{(1)}$ of   $\mathbf{u}$ is found to be
\begin{equation}
\label{Pr1}
\mathbf{u}^{(1)}=\lambda x \frac{\partial}{\partial x}+\gamma y \frac{\partial}{\partial y}+(\gamma-\lambda)y_x \frac{\partial}{\partial y_x}. 
\end{equation}
The vector field (\ref{Pr1}) represents an infinitesimal action of a one-parameter Lie group of transformations in a three-dimensional (prolonged) space. Hence, we expect to obtain $3-1 = 2$ fundamental differential invariants. Indeed, solving the corresponding partial differential equation by the method of characteristics, we arrive at the following set of two fundamenal differential invariants: 
\begin{equation}
\label{finvariants}
I_1 = \displaystyle  yx^{-\frac{\gamma}{\lambda}},\quad I_2=\displaystyle y_x x^{\frac{\lambda-\gamma}{\lambda}},
\end{equation}
as expected, which means that any other differential invariant of the prolonged  infinitesimal group action defined by (\ref{Pr1}) if a function of $I_1$ and $I_2$. Now, combining the fundamental differential invariants (\ref{finvariants}) in such a way that the parameters $\lambda$ and $\gamma$ disappear, we arrive at the following differential invariant: 
\begin{equation}
\label{ws}
{\bI}(I_1,I_2)=\frac{x y_x}{y},
\end{equation}
which we immediately recognize to be precisely the wage share $s_L$ (see, for example,  Rabbani \cite{SR2017} and Schneider \cite{DS2011} for more details). 

Therefore we conclude that not only the Cobb-Douglas production function (\ref{CD}), but also the wage share $s_L = {\bI}$ given by (\ref{ws}) is  a consequence of the exponential growth in $K$ and $L$ as a differential invariant obtained within the framework of the growth model $(G, \mathbb{R}^2_+)$, where the action of the Lie group  $G$  is given by (\ref{TPe}), that is 

$$\mbox{exponential growth} \Rightarrow \mbox{the wage share function (\ref{ws}).} $$

Now let us redo the above calculations  for the growth model $(G_1, \mathbb{R}^2_+)$, where the action of $G_1$ is given by (\ref{G1}) and thus give a solution to the seemingly unresolved problem of the determination of why  Bowley's law \cite{ALB1900, ALB1937} does not hold true anymore in post-1960s data \cite{SBSPG2003, MWLEBHAS2013, AG2007, LKBN2014}. 

First, we observe in the example considered above the exponential growth in $K$ and $L$ induced the corresponding exponential growth in the projective coordinates $x = L/K$ and $y = Y/K$. However,   the logistic growth in $K$ and $L$ given by (\ref{G1}) does not translate into the same type of transformations for the projective coordinates $x$ and $y$. Therefore, let us assume that the growth in $K$ is suppressed by, say,  excessive debt and so it does not affect logistic growth in $L$ and $Y$. Hence, both projective coordinates  $x$ and $y$ grow logistically, that is we have 
\begin{equation}
\label{logistic1}
\displaystyle \bar{x}=\dfrac{1}{1+(\frac{1}{x}-1)e^{-\lambda t}}, \quad \bar{y}= \dfrac{1}{1+(\frac{1}{y}-1)e^{-\gamma t}},  \,\, \lambda,  \gamma  \geqslant 0,
\end{equation}
 where we assumed without loss of generality that both carrying capacities  were equal to one. The corresponding infinitesimal action of the Lie group $G_1$ is given by the vector field
\begin{equation}
\label{u1}
\mathbf{u}_1 = \lambda x(1-x) \frac{\partial}{\partial x}+ \gamma y(1-y) \frac{\partial}{\partial y}.
\end{equation}
To determine its first prolongation $\mathbf{u}^{(1)}_1=\mbox{pr}^{(1)}\mathbf{u}_1$  we proceed as above  within the same framework as in the previous case (see the commutative diagram (\ref{diagram})).  We note first that the vector field $\mathbf{u}^{(1)}_1$ on $J^1 \pi$ is  projectable,  since the bundle $(T\mathbb{R}^2,\tau, \mathbb{R}^2)$ is endowed with a vector structure (see Saunders \cite{DJS1989}, Chapter 2 for more details). Next, define 
\begin{equation}\label{eq8}
\mathbf{u}^{(1)}_1=\lambda x(1-x)\frac{\partial}{\partial x} + \gamma y(1-y) \frac{\partial}{\partial y}+\xi(x,y,y_x) \frac{\partial }{\partial y_x}
\end{equation}
and require that the vector field (\ref{eq8}) is  a symmetry of the Cartan distribution, which will assure that (\ref{eq8}) is the first prolongation of (\ref{u1}). Indeed, consider again a basic contact form $\omega= \mbox{d}y-y_x \mbox{d}x$. Then again,  ${\bL}_{\mathbf{u}^{(1)}_1} (\omega)$ is a contact form iff $\mathbf{u}^{(1)}_1$ is a symmetry of the Cartan distribution on $J^1 \pi$, which in turn assures that (\ref{eq8}) is indeed the first prolongation of (\ref{u1}), where $\bL$ as before denotes the Lie derivative. Thus, we compute 

\begin{equation}
\label{com}
\begin{array}{ll}
{\bL}_{\mathbf{u}^{(1)}_1}(\omega) & = {\bL}_{\mathbf{u}^{(1)}_1} (\mbox{d}y-y_x\mbox{d}x) \\[0.3cm]
&= {\bL}_{\mathbf{u}^{(1)}_1}(\mbox{d}y)-({\bL}_{\mathbf{u}^{(1)}_1}(y_x)\mbox{d}x-y_x ({\bL}_{\mathbf{u}^{(1)}_1}(\mbox{d}x))\\ [0.3cm]
&=\mbox{d}(\mathbf{u}^{(1)}_1(y))-(\mathbf{u}^{(1)}_1(y_x))\mbox{d}x-y_x\mbox{d}(\mathbf{u}^{(1)}_1(x)) \\[0.3cm]
&= \gamma(1-2 y)\mbox{d} y-\xi(x,y,y_x)\mbox{d}x-\lambda(y_x\mbox{d}x-2xy_x\mbox{d}x) \\[0.3cm]
&=\gamma(1-2 y)(\omega+y_x\mbox{d}x)-(\xi(x,y,y_x)+\lambda y_x-2\lambda x y_x)\mbox{d}x \\[0.3cm]
&=\gamma (1-2 y)\omega+(\gamma y_x-2\gamma yy_x-\xi(x,y,y_x)-\lambda y_x+2 \lambda x y_x)\mbox{d}x. 
\end{array}
\end{equation}
In view of the above, ${\bL}_{\mathbf{u}^{(1)}_1}(\omega)$ is again a contact form, provided the expression in the parenthesis that appears in the last line of (\ref{com}) vanishes. Hence, we have
$$
\gamma y_x-2\gamma yy_x-\xi(x,y,y_x)-\lambda y_x+2 \lambda x y_x=0,
$$
or,
\begin{equation}
\label{xi}
\xi(x,y,y_x)=(\gamma-\lambda+2\lambda x-2\gamma y)y_x.
\end{equation}
We conclude therefore that the first prolongation of the vector field $\mathbf{u}_1$ given by (\ref{u1}) is the following fector field: 
\begin{equation}\label{eq7}
\mathbf{u}^{(1)}_1=\lambda x(1-x)\frac{\partial}{\partial x} + \gamma y(1-y) \frac{\partial}{\partial y}+(\gamma-\lambda+2\lambda x-2\gamma y)y_x \frac{\partial }{\partial y_x}, 
\end{equation}
whose infinitesimal action brings about the following two fundamental differential invariants: 
\begin{equation}
I_1=-\left(\frac{y-1}{y}\right) \left(\frac{x}{x-1}\right)^{\frac{\gamma}{\lambda}}, \quad
I_2=(2\gamma x)^2 \left(\frac{y_x}{(y-1)^2}\right) \left(\frac{1-x}{x}\right)^{\frac{\gamma+\lambda}{\lambda}}.
\end{equation} 
In order to eliminate the parameters $\lambda$ and $\gamma$ let us consider the following combination: 
\begin{equation}\label{com1}
{\bI}(I_1,I_2)=I_1 \cdot \frac{I_2}{(2\gamma)^2} =  x|x-1|  \left( \frac{y_x}{y |y-1|} \right). 
\end{equation}
\begin{definition}
The differential invariant $\bI$ given by (\ref{com1}) is  called a {\em modified wage share $s'_L = {\bI}$}, so that 
\begin{equation}
\label{ws1}
s'_L=\frac{|x-1|}{|y-1|}s_L=\mbox{const},
\end{equation}
where $s_L$ is the classical wage share given by (\ref{ws}). 
\end{definition}
\begin{remark}
\label{r7}
The modified wage share $s'_L$ given by (\ref{ws1}) is a differential invariant of the growth model $(G_1, \mathbb{R}_+^2)$, where the action of the Lie group $G_1$ is given by (\ref{G1}), while the classical wage share $s_L$ given by (\ref{ws}) {\em is not}. That is a reason why $s_L$ has been in decline: it may be attributed to the fact that  post-1960 economic data has been generated   within the framework of the growth model $(G_1, \mathbb{R}_+^2)$, rather than $(G, \mathbb{R}_+^2)$. More specifically, it follows that the decline in $s_L$ is due to the relation $\gamma > \lambda$ (see (\ref{ws1})). Indeed, if the output-to-capital ratio $y$ grows logistically faster than the labor-to-capital ratio $x$ under the condition of supressed capital (e.g.,  by excessive debt), that is if $\gamma >\lambda$ the ratio $\frac{|x-1|}{|y-1|}$ in (\ref{ws1}) clearly contributes to decline in $s_L$, since $s'_L$ is a constant. Simply put,  more wealth (real or perceived) distributed among fewer people implies a marked decrease in the classical wage share $s_L$ and so Bowley's law  \cite{ALB1900, ALB1937} no longer holds in the economic environment of the logistic growth model $(G_1, \mathbb{R}_+^2)$.  
\end{remark}

\begin{remark}
\label{r8}
The corresponding production function compatible with the infinitesimal action generated by the vector field $\mathbf{u}_1$ (\ref{u1}) is readily found to be 
\begin{equation}
\label{productionf1}
Y=f_{9}(K, L; t) = \frac{KL^{C_3}}{L^{C_3}+C_4 |L-K|^{C_3}},\quad C_3 \in (0,1), C_4 \in \mathbb{R},
\end{equation} 
which we derived by integrating the equation ${\bI} = \mbox{const}$, where ${\bI}$ is given by (\ref{com1}) and rewriting the solution in terms of $K$ and $L$.  

Now,  let us analyse the second new production function (\ref{productionf1}). The partial derivatives of the production function $f_9$ (\ref{productionf1}),  called in economic literature \emph{marginal productivities},  are found to be
\begin{equation}\label{MPK}
MP_K=\frac{1}{1+C_4 |1-\frac{K}{L}|^{C_3}}+C_3 C_4\frac{K}{L-K} \frac{|1-\frac{K}{L}|^{C_3}}{(1+C_4|1-\frac{K}{L}|^{C_3})^2}, 
\end{equation}
\begin{equation}\label{MPL}
MP_L=C_3C_4 \frac{K^2}{L(L-K)}\frac{|1-\frac{K}{L}|^{C_3}}{(1+C_4|1-\frac{K}{L}|^{C_3})^2}.
\end{equation}
Next, the slope of an isoquant is the \emph{marginal rate of technical substitution} (MRTS), or \emph{technical rate of substitution} (TRS). Thus,  $MRTS = \frac{MP_K}{MP_L}$ so that in our case 
\begin{equation}\label{MRTS1}
MRTS(K,L)=\frac{1}{C_3C_4}\frac{L(L-K)}{K^2}\frac{1+C_4|1-\frac{K}{L}|^{C_3}}{(1-\frac{K}{L})^{C_3}}+\frac{L}{K},
\end{equation}
which decreases when $L$ grows and $K$ declines.  We conclude, therefore, that \eqref{MRTS1} has concave up isoquants  when $L$ increases and $K$ decreases, that is  if the labour-capital ratio is less than approximately $\frac{1+C_3}{2}$, in which case  $MRTS$ increases,  while otherwise the isoquants are  concave down,  since $MRTS$ decreases.

Recall that the new productoin function (\ref{LPF1}) does not enjoy constant return to scale. Now let us examine the function (\ref{productionf1}) from this viewpoint. Indeed, for a factor $r>1,$ the substitution $(K, L) \rightarrow (rK,rL)$ in (\ref{productionf1}) yields
\begin{equation}
\begin{array}{ll}
f_9(rK,rL) & =\dfrac{rK(rL)^{C_3}}{(rL)^{C_3}+C_4|(rL)-(rK)|^{C_3}} \\
&=\dfrac{rKL^{C_3}}{L^{C_3}+C_4|L-K|^{C_3}}.
\end{array}
\end{equation}
which means that the new production function \eqref{productionf1} has \emph{constant returns to scale},  since it is a homogeneous function of degree one. Therefore we conclude that it satisfies {\em the law of diminishing marginal returns and has constant return to scale}, which means it has a great potential for playing a pivotal role  in various economic growth models. 

Finally, let us investigate the behavior of the new production function (\ref{productionf1}) as $t \to 0$ and $t\to\infty$ under the assumption that both  $K(t)$ and $L(t)$ grow logistically according to the one-parameter Lie group transformations defined by (\ref{G1}). To understand its behaviour when $K$ and $L$ are small, we employ economic reasoning. Thus, at the beginning of a production cycle a company, say, invests much of its resources into fixed assets (e.g., infrastructure, materials, land, etc) and so when $t$ is small it is safe to assume that $K \gg L$, which implies that 
\begin{equation}
f_9(t) \sim  \frac{1}{C_4} (K(t))^{1-C_3}(L(t))^{C_3},
\end{equation}
that is the production function $Y$ enjoys a similar behaviour to that of the Cobb-Douglas production function (\ref{CD}) that has constant returns to scale.  When $t \to \infty$ both $K$ and $L$ grow logistically and so we have by (\ref{productionf1})
$$\lim_{t\to\infty}f_9
(K, L; t) = \mbox{const}.$$ 

\end{remark}

\section{The new production function $f_5$ vis-\`a-vis economic data} 
\label{S7}

In this section we present a similar analysis to the one conducted by Cobb and Doublas \cite{CWCPHD1928}, namely we l compare the new production function with some available US economic data from 1947-2016.  We make use of the data from the period 1947-2016 that is provided by the Federal Reserve Bank of St. Louis (https://fred.stlouisfed.org),  employing the FRED tool.   The variables are as follows: $K$ --- capital services of nonfarm business sector \cite{K}, $L$ --- compensations of employees of nonfarm business sector \cite{L}, $Y$ --- real output of nonfarm business sector \cite{Y}. The values of all variables are dimensionless, they are index values with the values at  2009 taken as 100. To estimate the new production function (\ref{LPF1}), we have used R Programming \cite{R}, employing the method of least squares, and assuming the corresponding carrying capacities to be of the following values: $N_{f_5}= 120$, $N_L = 150$. We have also assumed that $\alpha + \beta = 1$. 

\begin{figure}[h]
 \centering
  \includegraphics[width=15cm,height=7cm]{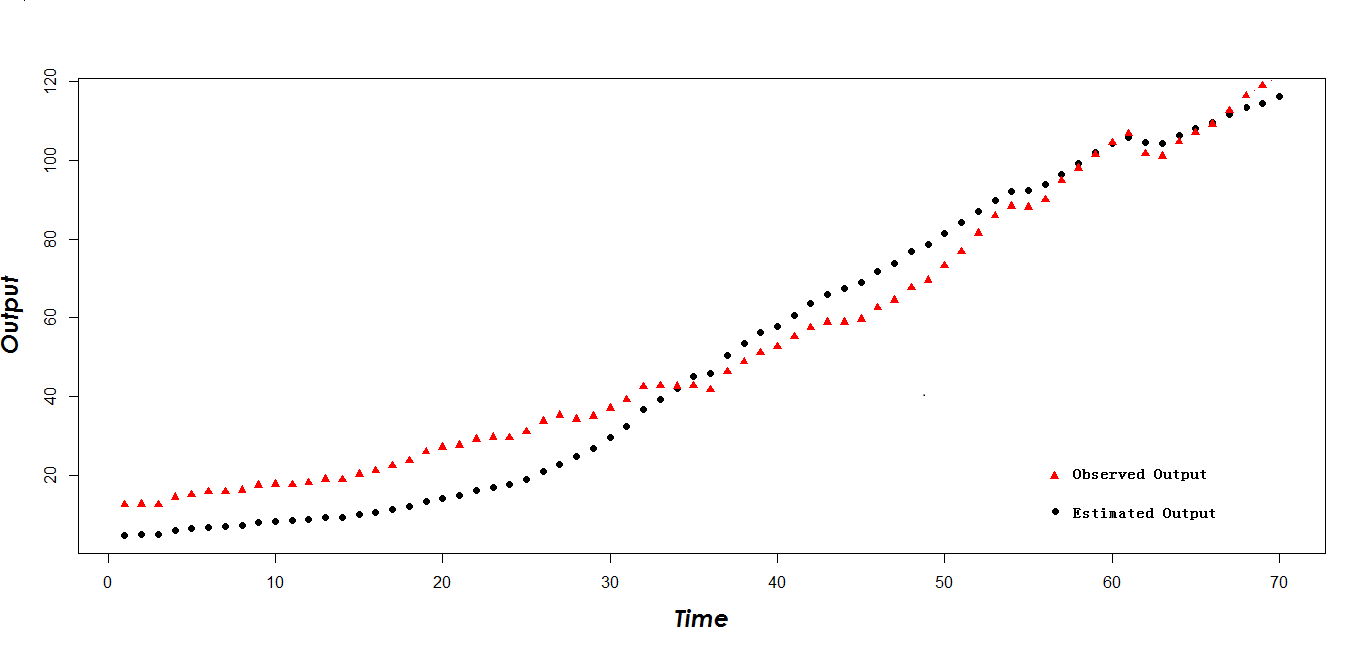}
   \caption{Observed output vs estimated output using  the new production function (\ref{LPF1}).} 
   \label{fige2}
\end{figure}

The resulting production function of the type (\ref{LPF1}) is found to be 
\begin{equation}
Y=\frac{120K^{(0.4063544)}L^{(0.5936456)}}{(0.3118901)|150-K|^{(0.4063544)}|150-L|^{(0.5936456)}+K^{(0.4063544)}L^{(0.5936456)}},
\end{equation}
where $C=0.3118901,$ $\alpha=0.4063544$ and $\beta=0.5936456$ (see Figure \ref{fige2}).

The elasticity of substitution $\sigma_1$ (see Sato \cite{RS1970}) of the new production function (\ref{LPF1}) in this case assumes the following form: 
\begin{equation}
\sigma_{1}=  \dfrac{\dfrac{\dot{L}}{L}-\dfrac{\dot{K}}{K}}{\dfrac{\dot{L}}{L}-\dfrac{\dot{K}}{K}-\dfrac{\dot{K}}{K-1}-\dfrac{\dot{L}}{L-1}},
\end{equation}
where $K=\frac{N_K C_1}{C_1+(N_K-C_1)e^{-a t}}$, $L=\frac{N_LC_2}{C_2+(N_L-C_2)e^{-b t}}$,   while $C_1$ and $C_2$ are constants.  The vairable   $\sigma_1$,  giving the best  estimate when $C_1=0.203,$ $a=0.129,$ $C_2=0.432$ and $b=0.118$,   ranges approximately from $-0.0151724079$ to $0.4982041724$. 

\begin{figure}[h]
 \centering
  \includegraphics[width=13cm,height=5cm]{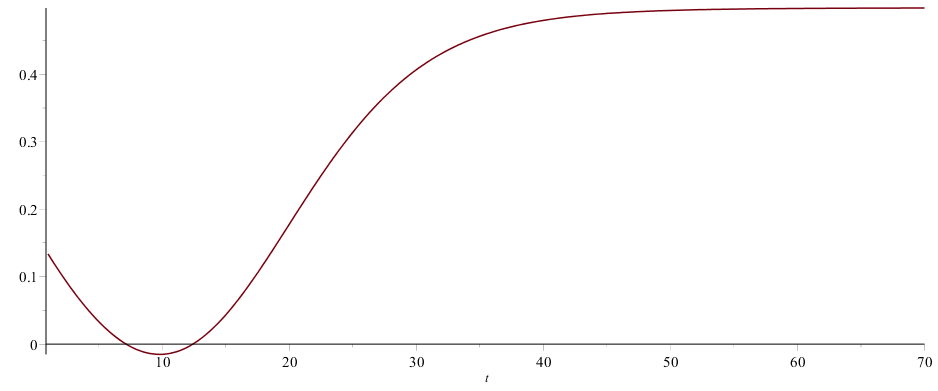}
   \caption{The elasticity of substitution of the new production function from 1947 to 2016.} 
   \label{fige3}
\end{figure}

Whether the function $f_5$,  derived using the Lie group theoretical methods,  can accurately  predict the future  still remains to be seen, but it looks like the function $f_5$  can ``predict'' the past. More specifically, while running our simulations, we have noticed that the negative value of  $\sigma_1=-0.0151724079$ occurs  in  the year of 1958 - excatly the year of a sharp economic downturn  \cite{RWG1959}, see Figure \ref{fige3}. 

\begin{figure}[h]
 \centering
  \includegraphics[width=15cm,height=7cm]{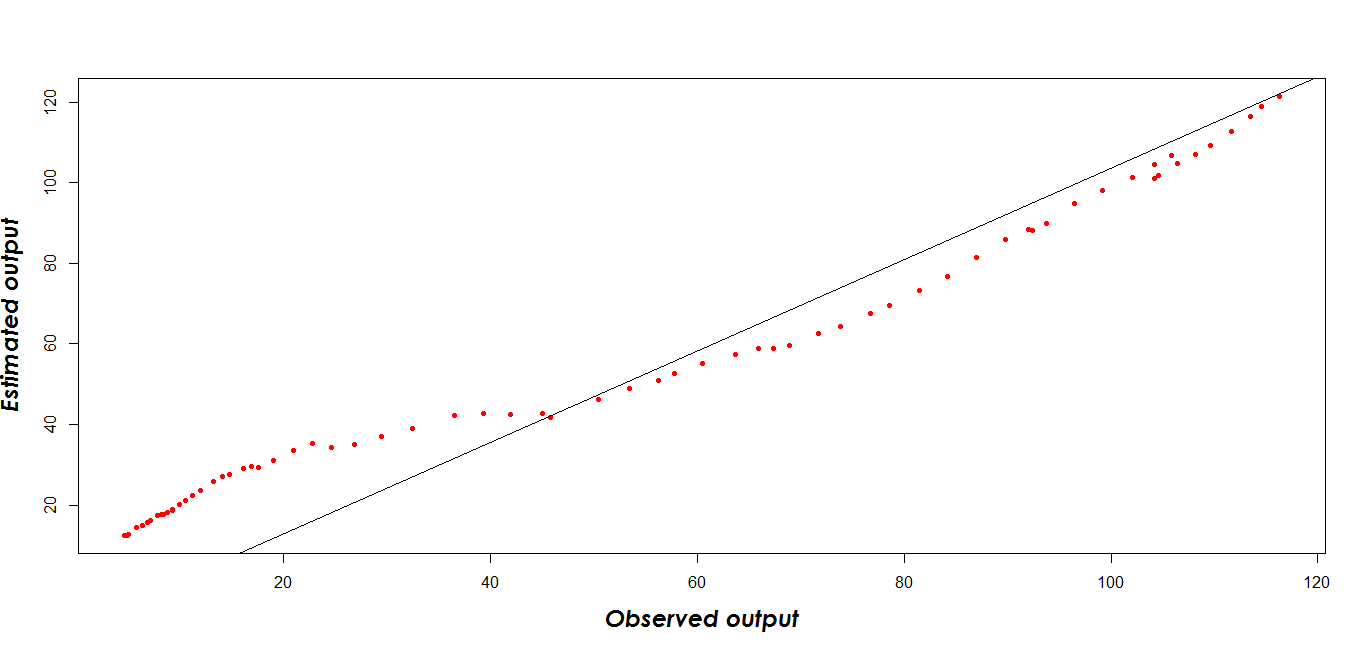}
   \caption{The linear regression of the observed and estimated outputs of the period from  1947 to 2016.} 
   \label{fige1}
\end{figure}

We conclude from the above that  the time series from the period 1947-2016 that compares the  observed and estimated outputs (see Figure \ref{fige1}) reveals that our model fits quite well the data with the the adjusted R-squared value of 97.65$\%$. On the other hand, the Cobb-Douglas function (\ref{CD}) with a constant elasticity of substitutions, i.e., $\sigma=1$, does not provide satisfactory results in terms of the values of parameters $\alpha$ and $\beta$. The best estimation of the Cobb-Douglas function that we managed to have obtained,  using the same method,  is as follows: 
\begin{equation}
Y=(0.2464455)K^{(1.6612365)}L^{(-0.6612365)},
\end{equation}
where $C=0.2464455,$ $\alpha=1.6612365$ and $\beta=-0.6612365$. We see that this (negative!) value of the parameter $\beta$ is not compatible with  the definition of the Cobb-Douglas production function given by  the formula (\ref{CD}). 

\section{Summary and discussion} 
\label{S8}

In this paper we have introduced a new (logistic)  growth model $(G_1, \mathbb{R}_+^2)$ given by (\ref{G1})  as an extension  and natural continuation of the preceeding studies in the area of economic growth done by Ryuzo Sato \cite{RS1970, RS1977, RS1980, RS1981}, as well as  a new framework for the development of more general production functions that we believe fit better  current economic data. The resulting new production functions (\ref{LPF1}) and (\ref{productionf1}) are consequences of the logistic growth in factors (capital and labor). The former function has shown to provide an adequate estimate for economic data, as for the latter --- there are indications that it will perform even better, the work in this direction is underway. Furthermore, we have presented a purely mathematical justification of why Bowley's law  \cite{ALB1900, ALB1937}  no longer holds true in post-1960 economic data by introducing a  new notion of  modified wage share  (\ref{ws1}). 

Our research has also demonstrated that {\em there can not be exponential growth of production while factors grow logistically}. We are inclined to believe that this is the most important consequence of our studies. Indeed, if one ``forces'' the production function  to grow exponentially  (i.e., by setting $H(f) = cf$   in (\ref{Covlogistic1})), while the factors $K$ and $L$ grow logistically as in (\ref{G1}),  the resulting production function will be of the form 
\begin{equation}
Y = f_{10}(K,L;t) = C_1\left(\frac{K}{|1-K|}\right)^{C_2} \left(\frac{L}{|1-L|}\right)^{C_3},
\label{productionf2}
\end{equation} 
where we assumed without loss of generality that $N_K = N_L = 1$. The production function $f_{10}$ (\ref{productionf2}) blows up very quickly near the singularities at $K=1$ and $L=1$. Similarly unsatisfactory result can by obtained by enforcing logistic growth in the production function, while the factors $K$ and $L$ grow exponentially, that is by setting $H(f) = cf(1-f)$ in (\ref{Cov1}): the resulting production function will not even grow. 

When we were starting this project, our original goal was to  only extend the theoretical framework based on the Lie group theory developed by Sato, we did not excpect that  the resulting production functions would perfom so well. Therefore the results obtained in this paper have exceeded our expectations. 

We see many  applicatoins in both economic theory of growth and applied mathematics  where the new production functions (\ref{LPF1}) and  (\ref{productionf1}), as well as the new modified wage share (\ref{ws1})  can be  used essentially  {\em mutatis mutandis} by simply replacing the Cobb-Douglas function or its generalizations (like the CES function, for example) and wage share with them as appropriate. 
 
As we have already mentioned in Introduction, the idea that exponential growth ought to be replaced with the logistic one is slowly but surely becoming more and more accepted by the scientists developing various growth models (see  Capasso {\em et al} \cite{CED2012}, Engbers {\em et al} \cite{EBC2014}, Brass \cite{WB1974},   Ferrara and Guerrini \cite{FG2008, FG2008a, FG2008b, FG2009, FG2009a, MFLG2009},  Leach \cite{DL1981}, Oliver \cite{ERO1982}, Tinter \cite{GT1952})  fore more details and references). 

 In light of the results that we have obtained so far,  some of the projects that we have learned from and appreciated so much, we belive could be modified accordingly,  which in turn may lead to more accurate mathematical models. For example, in Ferrara and  Guerrini \cite{FG2009a} the authors generalized the Ramsey model by introducing the logistic growth in $L$, which was  a very adequate assumption. However, they still used the Cobb-Douglas function which, we believe, is not entirely accurate, because the logistic growth in $L$ suggests that the growth model $(G_3, \mathbb{R}_+^2)$ given by (\ref{G3}) is underpinning the dynamics of the variables involved and so one has to use the corresponding production function compatible with (\ref{G3}) that is the function (\ref{Y3}) instead of the Cobb-Douglas production funciton (\ref{CD}). Similarly,  Capasso {\em et al} \cite{CED2012} did introduce a modified production function (\ref{LPF}) instead of the usual Cobb-Douglas production function (\ref{CD}), however it was done heuristically and a more natural choice for a production function in the model developed by the authors is either the production function (\ref{LPF1}) or (\ref{productionf1}), both of which were derived here in a systematic way. More specifically, the partial differential equation
$$\frac{\partial K}{\partial t}(x,t) = \Delta K(x,t)  + F(K(x,t), L(x,t)) - \delta K(x,t), \, x \in \Omega \subset \mathbb{R}^n, \, t \geqslant 0  $$
governing the dynamics of $K$ should use either (\ref{LPF1}) or (\ref{productionf1}) in place of $F$, which we believe will lead to more accurate results.

\section*{Acknowledgement}

The authors wish to thank Ryuzo Sato for valuable comments, constructive critique and suggestions. The second author (KW) wishes to thank Chaoyue Liu for his invaluable  help with R Programming  and  statistical analysis he used while working on the material presented in Section \ref{S7}.

\end{document}